\documentclass{article}

\usepackage{times}
\usepackage{amsmath}
\usepackage{amsthm}
\usepackage{amsfonts}
\usepackage{mathrsfs}
\usepackage{fullpage}

\title
{Interpolation and Sampling for Generalized Bergman spaces on finite
Riemann Surfaces}
\author
{Alexander P. Schuster and Dror Varolin}
\date{}

\newcommand{\scb}{{\mathscr B}}
\newcommand{\scf}{{\mathscr F}}
\newcommand{\sco}{{\mathscr O}}

\newcommand{\ms}{\medskip}

\newcommand{\dfn}[1]{{\noi {\bf Definition.\ } {\it #1}}}

\newcommand{\noi}{\noindent}
\newcommand{\pf}{{\noi \it Proof. }}
\newcommand{\rmk}{\noi {\bf Remark. }}

\newcommand{\ve}{\varepsilon}

\newcommand{\cb}{{\mathcal B}}
\newcommand{\cc}{{\mathcal C}}

\newcommand{\cv}{{\mathcal V}}

\newcommand{\fb}{{\mathfrak b}}

\newcommand{\ff}{{\mathfrak f}}

\newcommand{\vp}{\varphi}
\newcommand{\ii}{\sqrt{-1}}

\newcommand{\Z}{{\mathbb Z}}

\newcommand{\R}{{\mathbb R}}

\newcommand{\C}{{\mathbb C}}
\newcommand{\D}{{\mathbb D}}

\newcommand{\di}{\partial}
\newcommand{\dbar}{\overline \di}
\newcommand{\nbar}{\overline \nabla}

\newcommand{\relcomp}{\subset \subset}

\newcommand{\hs}{*}

\newcommand{\re}{{\rm Re}}
\newcommand{\im}{{\rm Im}}

\newcommand{\tensor}{\otimes}

\begin{document}

\maketitle

\newtheorem{thm}{Theorem}[section]
\newtheorem{mthm}{Theorem}
\newtheorem{lem}[thm]{Lemma}
\newtheorem{prop}[thm]{Proposition}
\newtheorem{cor}[thm]{Corollary}
\newtheorem{conj}[thm]{Conjecture}
\newtheorem{defn}[thm]{Definition}
\newtheorem{exa}[thm]{Example}

\section{Introduction}\label{intro}

In \cite{s1,sw,s2}  Seip {\it et al.}  characterized the sampling
and interpolating sequences in the Bargmann-Fock space of entire
functions that are square integrable with respect to the weight
function $e^{-|z|^2}$, and in the Bergman space of square integrable
holomorphic functions on the unit disk $\D$.

In both cases, the results are given in terms of densities. In the
complex plane $\C$, let $\Gamma$ be a discrete, uniformly separated
sequence, and define
\begin{eqnarray*}
D^+(\Gamma) := \limsup _{r \to \infty} \sup _{z\in \C} \frac{\#
(\Gamma \cap D(z,r))}{r^2} &\ & D^-(\Gamma) := \liminf _{r \to
\infty} \inf _{z\in \C} \frac{\# (\Gamma \cap D(z,r))} {r^2}.
\end{eqnarray*}
In the unit disk $\D$, the corresponding densities are defined in an
analogous manner, which is nevertheless slightly different.  Let
$\Gamma$ be a sequence that is uniformly separated in the
pseudo-hyperbolic distance. We then set
$$
D^+(\Gamma) := \limsup _{r \to 1} \sup _{z\in \D}
\frac{2\int_0^r\#(\Gamma\cap D(z,s))ds}{\log\frac{1}{1-r}}
$$
and
$$
D^-(\Gamma) := \liminf _{r \to 1} \inf _{z\in \D}
\frac{2\int_0^r\#(\Gamma\cap D(z,s))ds}{\log\frac{1}{1-r}},
$$ where $D(z,s)$ is the pseudohyperbolic disk of center $z$
and radius $s$. The numbers $D^{\pm}(\Gamma)$ are often called the
upper and lower densities of $\Gamma$. Seip {\it et al.} proved the
following, now celebrated theorem.

\medskip

\noi {\bf Theorem:} {\it A uniformly separated sequence $\Gamma$ is
an interpolating sequence for the Bargmann-Fock space or the Bergman
space if and only if $D^+(\Gamma) < 1$. It is a sampling sequence if
and only if $D^-(\Gamma) > 1$.}

\medskip

\noi The goal of the present article is to generalize the
sufficiency part of the Theorem of Seip {\it et al.} to the case of
open Riemann surfaces.  We fall somewhat short of this goal, but it
is not clear how short.  Indeed, we establish generalizations for
the case of finite Riemann surfaces.  However, the methods used do
generalize to other Riemann surfaces, and may even generalize to all
Riemann surfaces; we were unable to decide.

We now present our main results.  To this end, every open Riemann
surface admits a metric, locally denoted $e^{-2\nu} |dz|^2$, that we
call the {\it fundamental metric}. (See definition
\ref{fund-metric}.) Moreover, if the Riemann surface is hyperbolic,
then this metric is unique.   With the fundamental metric at hand,
we can associate to a smooth function $\vp : X \to \R$ and a
discrete subset $\Gamma \subset X$ two Hilbert spaces
$$
\scb ^2 = \scb ^2 (X,\vp) := \left \{ h \in \sco (X) \ ;\ ||h||^2 :=
\int _X |h|^2 e^{-2\vp} dA_{\nu} < +\infty \right \}
$$
and
$$
\fb ^2 = \fb ^2 (\Gamma ,\vp) := \left \{ (s_{\gamma})_{\gamma \in
\Gamma}\ ;\ \sum _{\gamma \in \Gamma} |s_{\gamma}|^2
e^{-2\vp(\gamma)} < +\infty\right \}.
$$

\dfn{A discrete set $\Gamma$ is said to be
\begin{enumerate}
\item an interpolation set if for every $(s_{\gamma}) \in \fb ^2$ there
exists $F \in \scb ^2$ such that for all $\gamma \in \Gamma$,
$F(\gamma) = s_{\gamma}$, and

\item a sampling set if there is a constant $M$ such that for all $F
\in \scb ^2$,
$$
\frac{1}{M} ||F||^2 \le \sum _{\gamma \in \Gamma} |F(\gamma)|^2
e^{-2\vp (\gamma)} \le M ||F||^2.
$$
\end{enumerate}
}

\noi For each locally integrable function $f : [0,R_X) \to
[0,\infty)$ and each $r \in (0,R_X)$, let $c_r := 2\pi \int _0 ^r
tf(t)dt$ and
$$
\xi _r (z,\zeta ) = \frac{1}{c_r} f(\rho _z(\zeta)) e^{2\nu} |d\rho
_z (\zeta)|^2 {\bf 1}_{D_r(z)}(\zeta),
$$
where ${\bf 1}_A$ denotes the characteristic function of a set $A$.
To every uniformly separated sequence (see section \ref{uss} for the
definition) we associate the upper and lower densities
$$
D^+_f (\Gamma) := \limsup _{r \to R_X} \sup _{z \in X} \sum _{\gamma
\in \Gamma} \frac{\frac{\pi}{2} \xi _r (\gamma ,z)}{e^{2\nu
(z)}\Delta \vp (z)}
$$
and
$$
D^-_f (\Gamma) := \liminf _{r \to R_X} \inf _{z \in X} \sum _{\gamma
\in \Gamma} \frac{\frac{\pi}{2} \xi _r (\gamma ,z)}{e^{2\nu
(z)}\Delta \vp (z)}
$$
where $\Delta$ is the Laplace operator.  Our main theorem can now be
stated as follows.
\begin{mthm}\label{main-thm}
Let $X$ be a finite Riemann surface, $\Gamma \subset X$ a uniformly
separated sequence, and $\vp : X \to R$ a subharmonic function such
that for some constant $C$, $\frac{1}{C} \le e^{2\nu}\Delta \vp \le
C$.  If $D^+_f(\Gamma) < 1$ the $\Gamma$ is an interpolation set,
while if $D^-_f (\Gamma) > 1$, then $\Gamma$ is a sampling set.
\end{mthm}

Partial results covering our theorem have been proved by others. For
the case of the plane but with more general subharmonic weights,
Theorem \ref{main-thm} was proved by Berndtsson-Ortega Cerd\`a
\cite{quimbo}.

Ohsawa \cite{o} has proved results on interpolation only, and in the
much more general context extension of $L^2$ holomorphic functions
from closed submanifolds of Stein manifolds.  His approach is
somewhat different than ours; he uses a method, pioneered by himself
and Takegoshi, of using a twisted $\dbar$ theorem at (an earlier
stage of the construction) to do the extension directly, rather than
use $\dbar$, as is done in \cite{quimbo}.  (Ohsawa argued later
\cite{o2} that his approach is more conducive to generalization than
the $\dbar$ approach used here.  We believe the methods of this
paper show that both approaches are equally generalizable.)

In the same paper \cite{quimbo}, Berndtsson and Ortega Cerd\`a also
treat the case of functions that are square integrable with respect
to a subharmonic weight on the unit disk.  A more careful analysis
shows that because of the curvature of the fundamental metric,
Theorem 1 does not cover all the cases those authors treat. A way to
compensate is to prove a second theorem in the case of the disk,
which allows some relaxation of the condition on $\Delta \vp$. This
is indeed what was done in \cite{quimbo}. Aesthetically, this
approach has the disadvantage of making the results of the disk and
the plane appear distinct.  Instead, in section \ref{metric-case} we
use conformal metrics on the Riemann surface to obtain more general
norms on our generalized Bergman spaces with metrics.  We then prove
results (Theorems \ref{interp-thm} and \ref{sampling-thm}) which
encompass Theorem \ref{main-thm}.

One additional novelty in our work is the introduction of a family
of densities, parameterized by locally integrable functions $f$ on
the positive real line.  This feature of our densities is likely to
be useful in applications, and gives new results even in the
classical cases of the plane and the disk.  We demonstrate some
useful consequences in the short Section \ref{examples-section} at
the end of the paper.

It is worth mentioning two additional things.

\begin{enumerate}
\item[(i)] As of right now we have not addressed
the question of necessity.  While we know that some of our density
conditions are not necessary for interpolation or sampling, we do
not know which, if any, are necessary.

\item[(ii)] Though we prove our theorems for finite Riemann surfaces,
our work applies to a much broader class of open Riemann surfaces.
In fact, we know of no example of a Riemann surface where our methods 
cannot be used to prove the corresponding version of Theorem \ref{main-thm}.

\end{enumerate}

The organization of our paper is as follows.  In section
\ref{finite-rs} we discuss potential theory and the resulting
analytic geometry of finite Riemann surfaces.  In section
\ref{metric-case} we introduce metrics into our scheme, and state
Theorems \ref{interp-thm} and \ref{sampling-thm}, which are the main
results of paper.  The hypotheses in those theorems appear rather
rigid, and thus it is not clear if or when they are satisfied.  Thus
in the same section we show that in fact there is always a
non-trivial case in which the hypotheses are satisfied.  In section
\ref{dbar} we discuss a $\dbar$ theorem, due to Ohsawa \cite{o2},
that will be used in the proof of the interpolation theorem
\ref{interp-thm}.  Since Ohsawa's theorem is more general, we give a
short, {\it ad hoc}  proof of the case we need here.  In section
\ref{cpct} we take a brief detour and establish interpolation and
sampling results for compact Riemann surfaces.  There are no
non-constant holomorphic functions on compact Riemann surfaces, so
we must look at sections of line bundles.  These spaces are always
finite dimensional, and sheaf theoretic methods give a complete
answer to the interpolation and sampling question.  Nevertheless we
prove a special case of the interpolation theorem in this setting,
to demonstrate how the $\dbar$ theorem is later used.  In section
\ref{proofs} we prove Theorems \ref{interp-thm} and
\ref{sampling-thm}.  Finally, in section \ref{examples-section} we
present a collection of examples of our main results in some special
cases.

\section{Analytic Geometry of Finite Riemann
surfaces}\label{finite-rs}

\subsection{Potential theoretic preliminaries}

We recall some basic, well known facts about fundamental solutions
of the Laplacian and about harmonic functions on Riemann surfaces.

\subsubsection{Extremal fundamental solutions and the fundamental metric}
\label{fund-sol-section}

We write
$$
\Delta := \ii \di \dbar
$$
for the Laplace operator.  Note that this is the complex analytic
convention, which is $1/4$ of the usual Laplace operator one
encounters in electromagnetism.  Let $\delta _{\zeta}$ denote the
Dirac mass at $\zeta$. The following definition is standard.

\ms

\dfn{The Green's function on a Riemann surface $X$ is the function
$G:X \times X \to [-\infty, 0)$ with the following properties.
\begin{enumerate}
\item[{\rm (a)}] For each $\zeta \in X$, $\Delta  _z G(z,\zeta) =
\frac{\pi}{2} \delta _{\zeta}(z).$ \item[{\rm (b)}] If $H : X \times
X \to [-\infty, 0)$ is another function with property {\rm (a)},
then $H(z,w) \le G(z,w)$ whenever $z\neq w$.
\end{enumerate}
}

\ms

\noi It can easily be deduced that the Green's function is
symmetric.

\medskip

Recall that a Riemann surface is said to be {\it hyperbolic} if it
admits a bounded subharmonic function, {\it elliptic} if it is
compact and {\it parabolic} otherwise.  It is well known that a
Riemann surface has a Green's function if and only if it is
hyperbolic.  Property (b) guarantees that the Green's function is
unique.

On the other hand, a Riemann surface admits an Evans kernel if and
only if it is parabolic (see page 352 of \cite{sario}).  Moreover,
after prescribing (with somewhat limited possibility) the
logarithmic singularity at infinity, the Evans kernel is unique up
to an additive constant.

\ms

\dfn{An Evans kernel on a Riemann surface $X$ is a symmetric
function $S:X\times X \to [-\infty ,+\infty)$ with the following
properties.
\begin{enumerate}
\item[{\rm (a)}] For each $\zeta \in X$, $\Delta _z  S(z,\zeta) =
\frac{\pi}{2} \delta _{\zeta}(z).$

\item[{\rm (b)}] For each $r \in \R$ and $p\in X$, the level set $\{
\zeta \in X\ ;\ S(\zeta,p) = r \}$ is compact and non-empty.
\end{enumerate}
}

\ms

We shall use the notation $E : X \times X \to [-\infty, R_X)$ to
denote either the Green's function or some chosen Evans kernel,
depending on whether the Riemann surface is hyperbolic or parabolic,
respectively. Define
$$R_X := \left \{
\begin{array}{c@{\quad}c}
1 & X \ {\rm is\ hyperbolic}\\
+\infty & X \ {\rm is\ parabolic}
\end{array} \right .
$$

Using the extremal fundamental solution we next define a notion of
distance on our Riemann surface, a distance that in general fails to
satisfy the triangle inequality.

\begin{defn}\label{fund-metric}
Let
$$
\rho _z (\zeta) := e^{E(z,\zeta)}, \quad D_{\ve} (z) := \left \{
\zeta \in X\ ;\ \rho _z (\zeta) < \ve \right \} \quad {\rm and}
\quad S_{\ve} (z) = \di D_{\ve} (z) .
$$
The fundamental metric $e^{-2\nu}$ is given by the formula
$$
e^{-2\nu(z)}|dz|^2 = \lim _{\zeta \to z} |\di \rho _z(\zeta)|^2.
$$
\end{defn}

\subsubsection{Green's Formula and mean values}

Recall that on a Riemann surface with a conformal metric, the Hodge
star operator simplifies somewhat when expressed in analytic
coordinates $z=x+\ii y$:  if $f$ is a real-valued function, $\alpha
= \alpha _1 dx + \alpha _2 dy$ is a real 1-form and $\vp dx\wedge
dy$ is a real 2-form, then one has
\begin{eqnarray*}
*f &=& f dA_g = e^{-2\psi} f dx \wedge dy \\
*\alpha &=& -\alpha _2 dx + \alpha _1 dy \\
*(\vp dx \wedge dy )&=& e^{2\psi} \vp.
\end{eqnarray*}
Using this, we have $4\Delta = d*d$ (recall that $\Delta = \di \dbar
= \frac{1}{4}(\di _x^2 + \di _y ^2)$ in our convention), and Green's
formula can be written
\begin{eqnarray}\label{grn-fa}
4 \int _Df\Delta h-h\Delta f=\int _{\di D}f * dh - h * df.
\end{eqnarray}

Let $X$ be an open Riemann surface and $Y \relcomp X$ an open
connected subset whose boundary consists of finitely many smooth
Jordan curves.  It is well known that the Green's function $G_Y$ for
$Y$ exists and is continuous up to the boundary.  Moreover, the
exterior derivative $d(G_Y (\zeta , \cdot ))$ is also continuous up
to the boundary.

\medskip

\rmk {\it One can construct the Green's function $G_Y$ from the
extremal fundamental solution $E$ of $X$ as follows.  Since $Y$ has
smooth boundary, the Dirichlet Problem of harmonic extension from
the boundary can be solved on $Y$.  We then take
$$G_Y (\zeta ,z) := E(\zeta, z) - h_{\zeta}(z),$$
where $h_{\zeta}$ is the harmonic function in $Y$ that agrees with
$E(\zeta,\cdot)$ on the boundary of $Y$.}

\medskip

We write
$$H_{r,\zeta}(z) := G_{D_r(\zeta)}(\zeta , z).$$
In fact, the function $H_{r,\zeta}$ has a particularly simple form
in terms of the extremal fundamental solution $E$:
\begin{eqnarray}\label{h=e}
H_{r,\zeta}(z) = E(z,\zeta) - \log r, \quad z \in D_r (\zeta).
\end{eqnarray}
Moreover, in this case we don't need to assume that $r$ is a regular
value of $\rho _{\zeta}$.

Putting $D = D_r(z)$ and $h = H_{r,z}$ in (\ref{grn-fa}) and using
the definition of Green's function, we obtain the following lemma.
\begin{lem}\label{mvt}
Let $r < R_E$ and $\zeta \in X$. Then
\begin{eqnarray}\label{pois}
2\pi f(z) = \int _{S_r(z)} f \hs dE_z + \int _{D_r(z)} H_{r,z}
\Delta f.
\end{eqnarray}
In particular, if $f$ is subharmonic, then
\begin{eqnarray}\label{sub-pois}
f(z) &\le& \frac{1}{2\pi}\int _{S_r(z)} f \hs dE_z
\end{eqnarray}
with equality when $f$ is harmonic.
\end{lem}

\subsection{Finite Riemann surfaces}

\subsubsection{Definition and construction of finite Riemann
surfaces}

Recall that a finite Riemann surface is a two dimensional compact
manifold with boundary, possibly with a finite number of points
removed.  Thus the topological data determining the Riemann surface
is finite, hence the name.

There are two types of finite Riemann surfaces.  One type has only
punctures and no one dimensional components, while the other type
does have smooth boundary components.  The first type of is always
parabolic (unless it has no punctures, in which case it is elliptic)
while the second type is always hyperbolic.

An alternate description of a finite Riemann surface $X$ can be
given as follows:  $X$ is a (not necessarily compact) manifold with
compact boundary, and in addition $X$ can be decomposed as
$$
X = X_{\rm core} \cup \bigcup _{j=1} ^N U_j,
$$
where $X_{\rm core}$ is a compact manifold with smooth boundary, and
each $U_j$ is biholomorphic to a punctured disk whose outer boundary
is one of the smooth boundary curves of $X_{\rm core}$.  (Of course,
$X_{\rm core}$ may have some other boundary components that do not
meet one of the $U_j$.)  The $U_j$ correspond to the punctures.

While every finite Riemann surface with no one dimensional boundary
is obtained from a compact Riemann surface by removal of a finite
number of points, there is an almost equally simple way to construct
hyperbolic finite Riemann surfaces; simply take a compact Riemann
surface and remove a finite number of smooth Jordan curves so that
the resulting surface as two components. Then either component is a
finite Riemann surface, and one can further remove any finite number
of points.

In fact, all finite Riemann surfaces are of this type.  Indeed, we
can fill in the punctures complex analytically (since they are just
punctured disks) to obtain a compact Riemann surface with boundary
$$
\tilde X =X_{\rm core} \cup \bigcup _{j=1} ^N \overline{U_j},
$$
and then form the so-called double of $\tilde X$.  For more on this
well-known construction see, for example, \cite{ss}.

\subsubsection{Analytic-Geometric properties of finite Riemann surfaces}

We shall now derive certain analytic-geometric properties of finite
Riemann surfaces that are useful in the proofs of our main theorems.

\begin{thm}\label{finite-fund-finite}
Let $X$ be a finite Riemann surface with extremal fundamental
solution $E$.  Then for each sufficiently small $\sigma \in (0,R_X)$
there is a constant $C=C_{\sigma}$ such that for all $z \in X$ and
all $\zeta \in D_{\sigma} (z)$ the following estimate holds.
\begin{eqnarray}\label{ff-est}
\frac{1}{C} \le e^{2\nu} |\di \rho _z (\zeta)|^2 \le C.
\end{eqnarray}
\end{thm}

Before getting to the proof, we make a few observations. Suppose
that we are given an extremal fundamental solution $E$ on our
Riemann surface $X$.  By the definition of a fundamental solution of
the Laplacian, if $z$ is a local coordinate on $U \subset X$ then
there exists a function $h_U(\zeta,\eta)$, harmonic in each variable
separately, such that $h_U(\zeta ,\eta)=h_U(\eta,\zeta)$ and
$$
E(p,q) = \log |z(p) -z(q)| + h_U (z(p),z(q)).
$$
The dependence of $h_U$ on $z$ is determined by the fact that $E$ is
globally defined.

For simplicity of exposition, we abusively write
$$E(z,\zeta) = \log |\zeta -z| + h(z,\zeta).$$
We then have that $\rho _z(\zeta) = |z-\zeta|e^{h(z,\zeta)}$ and,
differentiating, we obtain
$$
\di \rho _z(\zeta) = \frac{\overline{\zeta - z}}{|\zeta
-z|}e^{h(z,\zeta)} \left ( \frac{1}{2} + (z-\zeta) \di _{\zeta}
h(z,\zeta) \right ).
$$
It follows that
\begin{eqnarray}\nonumber
e^{-2\nu (\zeta)} &=& e^{2h(\zeta,\zeta)}, \qquad {\rm and} \\
\label{h-relation} e^{2\nu (\zeta)}\left | \di \rho _z (\zeta)
\right |^2 &=& e^{2(h(z,\zeta) -h(\zeta,\zeta))}\left | 1 + 2 (\zeta
-z) \frac{\di h(z,\zeta)}{\di \zeta} \right |^2.
\end{eqnarray}
We point out that this implies in particular  that the right hand
side of (\ref{h-relation}) is well defined, since this is the case
for the left hand side.

\medskip

\noi {\it Proof of theorem \ref{finite-fund-finite}.} We shall break
up the proof into the hyperbolic and parabolic case.

\medskip

\noi {\it The case of bordered Riemann surfaces.}

We realize $X$ as an open subset of its double $Y$.  Since
$\overline X = X\cup \di X$ is compact, it suffices to bound the
right hand side of (\ref{h-relation}) in a set $U\cap X$, where $U$
is a coordinate chart in $Y$.  For coordinate charts whose closure
lies in the interior $X$, it is clear that this can be done. Indeed,
if $U \relcomp X$ and $z,\zeta \in U$, then $h$ is a smooth function
that is harmonic in each variable separately, and $\rho _z(\zeta)
\asymp |\zeta-z|$ uniformly on $U$.  Thus by taking $\sigma$
sufficiently small, we obtain the estimate (\ref{ff-est}) for all $z
\in U$ and $\zeta \in D_{\sigma}(z)$.  We thus restrict our
attention to the boundary.

There are two types of boundary points; zero dimensional and one
dimensional.  However, the Green's function ignores isolated zero
dimensional boundary components, since they have capacity zero.  (In
particular, the distance $\rho _z$ fails to be proper when there
are punctures.)  Thus we may assume that there are no punctures.

Let $U \subset Y$ be a coordinate neighborhood of a boundary point
$x \in \di X$.  By taking $U$ sufficiently small, we may assume that
$U$ is the unit disk in the plane, that $U\cap X$ lies in the upper
half plane and that $\di X$ lies on the real line.  It follows that
the Green's function is given by
$$
E(z,\zeta) = \log |z-\zeta| - \log |\bar z - \zeta| +F(z,\zeta),
$$
where $F(z,\zeta)$ is smooth and harmonic in each variable on a
large open set containing the closure of $U$.  Indeed, the Green's
function for the upper half plane is $\log |z-\zeta| - \log |\bar z
- \zeta|$.  The regularity of $F$ then follows from the construction
of Green's functions on finite Riemann surfaces using harmonic
differentials on the double. (See \cite{ss}, $\S$4.2.) It follows
that in $U$,
$$2\frac{\di h(z,\zeta)}{\di \zeta} = -\frac{1}{\bar z - \zeta}  +
2\frac{\di F (z, \zeta)}{\di \zeta} \qquad {\rm and} \qquad \rho _z
(\zeta) \ge C \frac{|z-\zeta|}{|\bar z - \zeta|}.$$ Thus
\begin{eqnarray*}
\left |2(\zeta -z)\frac{\di h(z,\zeta)}{\di \zeta} \right | &\le&
\frac{|z-\zeta|}{|\bar z- \zeta|} + 2 |z-\zeta| \left |
\frac{\di F(z,\zeta)}{\di \zeta} \right | \\
&\le& C \frac{|z-\zeta|}{|\bar z- \zeta|} \\
&\le & C' \rho _z (\zeta),
\end{eqnarray*}
where the constant $C'$ depends only on the neighborhood $U$. The
proof in the hyperbolic case is thus complete.

\medskip

\noi {\it The case of compact Riemann surfaces with punctures.} Let
$E$ be the Evans kernel of $X$.  Fix $p \in X$ and choose $r$ so
large that the set $X - D_r(z)$ is a union of punctured disks
$U_1,..., U_N$.  We may think of each $U_j$ as sitting in $\C$, with
the puncture at the origin.

Since $D_r(z) \relcomp X$, each $x \in \overline{D_r(z)}$ has a
neighborhood $U$ for which the expression (\ref{h-relation}) is
bounded above and below by positive constants, depending only on
$U$, whenever $\rho _z(\zeta) < \sigma$ for some sufficiently small
$\sigma$ again depending only on $U$. Indeed, in any such
neighborhood the function $h$ is very regular, and $\rho _z(\zeta)$
is uniformly comparable to $|z-\zeta|$.

For $z,\zeta \in U_j$, the Evans kernel has the form
\begin{eqnarray}\label{ek-local}
E(z,\zeta) = \log |z-\zeta| - \lambda _j \log |\zeta| + F(z,\zeta),
\end{eqnarray}
where $\lambda _j > 0$ with $\lambda _1 + ... + \lambda _N = 1$, and
$F(z,\zeta)$ is smooth across the origin (see \cite{sario}). Indeed,
using the method of constructing harmonic differentials with
prescribed singularities (see \cite{ss} $\S$2.7) we can construct a
function with the right singularities, defined everywhere on $\bar
X$.  Such a function clearly can be written in the form
(\ref{ek-local}) near the puncture.  Thus by the uniqueness of the
Evans kernel with prescribed singularities at the punctures, this
function must differ from $E$ by a constant.

It follows that in $U$,
$$2\frac{\di h(z,\zeta)}{\di \zeta} = -\frac{\lambda _j}{\zeta}  +
\frac{\di F (z, \zeta)}{\di \zeta}$$ and
$$\rho _z (\zeta) \ge C \frac{|z-\zeta|}{|\zeta|}.$$
Thus
\begin{eqnarray*}
\left |2(\zeta -z)\frac{\di h(z,\zeta)}{\di \zeta} \right | &\le&
\lambda _j \frac{|z-\zeta|}{|\zeta|} + 2 |z-\zeta| \left |
\frac{\di F(z,\zeta)}{\di \zeta} \right | \\
&\le& C \frac{|z-\zeta|}{|\zeta|} \\
&\le & C' \rho _z (\zeta),
\end{eqnarray*}
where again the constant $C'$ depends only on the neighborhood $U$.
The proof of Theorem \ref{finite-fund-finite} is thus complete.\qed

\begin{prop}\label{hhom}  Let $X$ be a finite Riemann surface.  Then there exists
a constant $C$ such that, for sufficiently small $\sigma > 0$ and
all $z \in X$,
\begin{eqnarray}\label{wp}
&& \sup _{w\in D_{\sigma}(z)}\exp \left ( \frac{4}{\pi} \int
_{D_{2\sigma}(z)}\!\!\!\!\!\!\!\!\!\!\!\!\! -G(w,\zeta)
e^{-2\nu(\zeta)}
\right ) \\
\nonumber && \qquad \qquad \qquad \le C \inf _{w\in
D_{\sigma}(z)}\exp \left ( \frac{4}{\pi} \int
_{D_{2\sigma}(z)}\!\!\!\!\!\!\!\!\!\!\!\!\! -G(w,\zeta)
e^{-2\nu(\zeta)} \right ) < +\infty,
\end{eqnarray}
where $G$ is the Green's function for the domain $D_{2\sigma}(z)$.
\end{prop}

\noi {\it Sketch of proof.} Once again we can use compactness
properties of finite surfaces.  The finiteness of the integrals in
question is easy, since extremal fundamental solutions have only a
logarithmic singularity, and are thus locally integrable. Thus we
restrict ourselves to estimating near the boundary.

The local analysis used in the proof of Theorem
\ref{finite-fund-finite} shows that, near the boundary, the disks
$D_{\sigma}(z)$ are simply connected and that the metric $e^{-2\nu}$
is equivalent to the Poincar\'e metric of the disk in the hyperbolic
case, and the metric $|z|^{-2} |dz|^2$ in the parabolic case.

The hyperbolic case follows from the fact that the Green's function
$G(w,\zeta)$ is comparable to the Green's function of the disk. In
the parabolic case it is easier to work with the complement of the
unit disk rather than the punctured disk.  Then the metric
$e^{-2\nu}$ is comparable to the Euclidean metric, the Green's
function $G(w,\zeta)$ is comparable to the Green's function of the
plane, and the necessary estimate follows as in the Euclidean case.
This completes the sketch of proof.  \qed

\medskip

We shall also have use for the following lemma.

\begin{lem}\label{hh}
Let $X$ be a finite Riemann surface.  Let $\sigma >0$ be a fixed,
sufficiently small constant. If $\vp$ is a function for which
$e^{2\nu} \Delta \vp$ is bounded above and below by positive
constants, then there is a constant $C=C_{\sigma}$ such that, for
all $z \in X$ and all $w \in D_{\sigma}(z)$,
\begin{eqnarray}\label{hh-est}
\exp \left ( \frac{4}{\pi} \int
_{D_{2\sigma}(z)}\!\!\!\!\!\!\!\!\!\!\!\!\! -G(w,\zeta) \Delta
\vp(\zeta) \right ) \le C
\end{eqnarray}
\end{lem}

\pf By Theorem \ref{finite-fund-finite}, Proposition \ref{hhom} and
the boundedness of $\Delta \vp$, it suffices to prove the result
when $\Delta \vp (\zeta)= |d\rho _z(\zeta)|^2$ and $w=z$. In this
case, it is easy to show that the integral is equal to $8\sigma ^2$.
\qed

\medskip

The next result we will need is a global version of the Cauchy
estimates on a Riemann surface with Riemannian metric.

\begin{prop}\label{global-cauchy-est}
Let $X$ be a finite Riemann surface and let $g$ be a conformal
metric for $X$.  Then for every $\sigma \in [0, R_X)$ and $\ve
> 0$ there exists a constant $C_{\ve ,\sigma}$ such that
for any $x \in X$ the following Cauchy estimates hold.
\begin{eqnarray}\label{est}
\sup _{D_{\ve} (x)} |h|^2 &\le& C_{\ve,\sigma} \int
_{D_{\sigma}(x)}\!\!\!\!\!\!\!\! |h|^2 dA_g,
\end{eqnarray}
and
\begin{eqnarray}\label{der-est}
\sup _{D_{\ve} (x)}|\di \rho _{x}|^{-2}|h'|^2 &\le & C_{\ve, \sigma}
\int _{D_{\sigma}(x)}\!\!\!\!\!\!\!\! |h|^2 dA_g.
\end{eqnarray}
\end{prop}

\pf  We begin with the following lemma.

\begin{lem}\label{kg-kernel}
Let $X$ be a finite Riemann surface.  Then for every $x \in X$ there
exists a function $K^x : X \times X \to \R$ such that the following
hold for any $\sigma \in [0,R_X)$:
\begin{enumerate}
\item In the sense of distributions, $\Delta _z K^x (z,\zeta) =
\frac{\pi}{2} \delta _z(\zeta) $ for all $z,\zeta \in
D_{\sigma}(x)$. \item For every $\ve < \sigma /4$ there exists a
constant $C_{\ve ,\sigma}$ such that for any $x \in X$ the following
estimates hold:
\begin{eqnarray}
\label{c1} &&\sup _{z\in D_{\ve} (x)} \int _{V_{\sigma}(x)}
e^{2\psi}\left | \frac{\di \rho _{x}}{\di \zeta} \frac{\di
K^x(z,\zeta)}{\di \zeta} \right | ^2 \le
C_{\ve,\sigma}\\
\label{c2} && \sup _{z\in D_{\ve} (x)}|\di \rho _x (z)|^{-2} \int
_{V_{\sigma}(x)} e^{2\psi}|\di \rho _x |^{2}\left | \frac{\di ^2
K^x(z,\zeta)}{\di z \di \zeta} \right | ^2 \le C_{\ve,\sigma}
\end{eqnarray}
\end{enumerate}
Here $V_{\sigma}(x) := D_{\sigma}(x) -D_{\sigma /2}(x)$.
\end{lem}

\noi {\it Sketch of proof.} In the case of a bordered Riemann
surface with a finite number of punctures, one can find a function
$K^x$ that does not depend on the point $x$.  This is done as
follows. Let $Y$ be the double of $X$, and fix any smooth distance
function on $Y$. We let $X_{\ve}$ be the set of all $x \in Y$ that
are a distance less than $\ve$ from $X$. For $\ve$ sufficiently
small, $X_{\ve} - X$ is a finite collection of annuli whose inner
boundaries form the boundary of $X$.  We may take for our
Cauchy-Green kernel the Green's function of $X_{\ve}$. We leave it
to the reader to check that the relevant estimates hold.

In the case of an $N$-punctured compact Riemann surface, one
decomposes $X$ as
$$X= X_{\rm core}\cup \bigcup _{j=1} ^N U_j,$$
where $X_{\rm core}$ is a bordered Riemann surface, and each $U_j$
is a neighborhood of a puncture biholomorphic to the punctured disk.
Each surface in the union has a Cauchy-Green kernel by the
construction in the bordered Riemann surface case, and thus we are
done.\qed

\medskip

\noi {\it Completion of the proof of Proposition
\ref{global-cauchy-est}.} Let $f \in \cc ^{\infty} _0
(D_{\sigma}(x))$ and write $K^x_z(\zeta) = K^x(z,\zeta)$.  Applying
formula (\ref{grn-fa}) with $h (\zeta) = K^x _z(\zeta)$, we obtain
$$
\frac{\pi}{2} f(z) = \int _{D_{\sigma}(x)}\!\!\!\!\!\!\!\!\! K^x_z
d\dbar f = \int _{D_{\sigma}(x)}\!\!\!\!\!\!\!\!\! \dbar f \wedge
\di K^x_z.
$$
Now let $\ve < \sigma / 4$ and let $\chi \in \cc ^{\infty}
_0([0,3\sigma /4))$ be such that
$$
\chi |[0,\sigma /2] \equiv 1 \quad {\rm and}\quad  \sup |\chi '| \le
\frac{5}{\sigma}.
$$
If $h \in \sco (D_{\sigma}(x))$, then with $z \in D_{\ve} (x)$ we
have
\begin{eqnarray}\label{pre-c}
h(z) = \int _{D_{\sigma}(x)} \!\!\!\!\!\!\!\!\! h \chi ' (\rho _{x})
\dbar \rho _{x} \wedge \di K^x_z.
\end{eqnarray}
An application of the Cauchy-Schwarz inequality and the estimate
(\ref{c1}) gives the inequality (\ref{est}), while differentiation
of (\ref{pre-c}) followed by an application of the Cauchy-Schwarz
inequality and the estimate (\ref{c2}) gives inequality
(\ref{der-est}).\qed

\medskip

\rmk {\it Note that were it not for the requirement that
$C_{\ve,\sigma}$ be independent of
 $x$, Proposition \ref{global-cauchy-est} would follow without
{\rm (\ref{c1})} and {\rm (\ref{c2})}.}

\subsubsection{Discrete subsets in finite Riemann
surfaces}\label{uss}

Let $X$ be an open Riemann surface. Our work on sampling and
interpolation sequences requires the notion of the separation of a
sequence.  For a measurable subset $A \subset X$, let
$$
D_r(A)=\{w\in X\ ;\ w\in D_r(a)\ {\rm for\ some\ } a\in A\}.
$$
We define two separation conditions on a sequence $\Gamma$, both of
which are given in terms of the distance induced by the extremal
fundamental solution.

\begin{defn}\label{usep-and-sparse} Let $\Gamma \subset
X$ be a discrete set.
\begin{enumerate}

\item The separation constant of $\Gamma$ is the number
\begin{eqnarray*}
\sigma (\Gamma) &:=& \sup \{ r\ ;\ D_r(\gamma) \cap D_r(\gamma ') =
\emptyset \},
\end{eqnarray*}
and say that $\Gamma$ is uniformly separated if $\sigma (\Gamma) >
0$.

\item We say $\Gamma$ is sparse if there is a positive constant
$N_{r,\ve}$, depending only on $0<r,\ve<R_X$, such that the number
of points of $\Gamma$ lying in the set $D_r(D_\ve(z))$ is at most
$N_{r,\ve}$ for all $z\in X$.
\end{enumerate}
\end{defn}
In both the complex plane and the unit disk, the triangle inequality
allows one to easily show that a uniformly separated sequence is
sparse.

In both of these situations, the triangle inequality allows one to
estimate the diameter of a set $D_{\ve}(D_r(a))$ in terms of $\ve$
and $r$, uniformly for all $a$.

Such an estimate can always be found if it is allowed to depend on
the base point $a$.  This situation can be made uniform when $X$ is
a finite Riemann surface. As in the proofs of Theorem
\ref{finite-fund-finite} and Propositions \ref{hhom} and
\ref{global-cauchy-est}, we can take advantage of the compactness in
the picture.  In particular, we have uniform estimates if we have
them in neighborhoods of the boundary.  But on the boundary, the
potential theory of $X$ is either like that (near the boundary) in
the upper half plane or (near infinity) in the plane, where we know,
from triangle inequalities in those cases, that the needed estimates
hold. We thus have the following proposition.

\begin{prop}\label{sep->sparse}
In a finite Riemann surface $X$ every uniformly separated sequence
is sparse.
\end{prop}

We do not know whether Proposition \ref{sep->sparse} holds if one
removes the finiteness condition.

\section{Bergman spaces with metrics}\label{metric-case}

It is also interesting to introduce, in addition to the weight in
question, a metric.  Thus, suppose in addition to an open Riemann
surface $X$, a discrete subset $\Gamma$ and a weight function $\vp$,
we are also given a conformal metric $g$.  Thus we modify our
Hilbert spaces as follows:
$$
\scb ^2 = \scb ^2 (X,g,\vp) := \left \{ h \in \sco (X)\ ;\ ||h||^2
:= \int _X |h|^2 e^{-2\vp} dA_g < +\infty \right \},
$$
and
$$
\fb ^2 = \fb ^2 (\Gamma ,g,\vp) := \left \{ (s_{\gamma})_{\gamma \in
\Gamma}\ ;\ \sum _{\gamma \in \Gamma} |s_{\gamma}|^2 e^{-2\vp
(\gamma)} A_g(D_{\sigma}(\gamma))<+\infty \right \},
$$
where $A_g(B)=\int _B dA_g$.

\ms

As before, we say that a uniformly separated sequence $\Gamma$ is
interpolating if for any $(s_{\gamma})\in \fb ^2$ there exists $F\in
\scb ^2$ such that for all $\gamma \in \Gamma$, $F(\gamma) =
s_{\gamma}$.  On the other hand, the sequence $\Gamma$ is sampling
if there exists a constant $M$ such that for all $F \in \scb ^2$,
\begin{eqnarray}\label{sampling-ineq}
\frac{1}{M} ||F||^2 \le \sum _{\gamma \in \Gamma} |F(\gamma)|^2
e^{-2\vp (\gamma)} A_g (D_{\sigma}(\gamma)) \le M ||F||^2.
\end{eqnarray}

To obtain sufficient conditions for interpolation and sampling
sequences, the definition of the densities must be changed slightly.

\subsection{The definition of the upper and lower densities} We
associate to our metric $g=e^{-2\psi}$ the two functions
$$
u_{\psi} := \psi - \nu \quad {\rm and} \quad \tau _{\psi} :=
e^{2\nu} \left (\Delta \psi + 2 \Delta u_{\psi} -2|\di u_{\psi}|^2
\right ).
$$
For each locally integrable function $f:[0,R_X) \to [0,\infty)$ and
each $r \in \left (0, R_X\right )$, we associate to every uniformly
separated sequence $\Gamma$ upper and lower densities, defined by
\begin{eqnarray}\label{d+}
D^+_f(\Gamma) := \limsup _{r \to R_X} \sup _{z\in X} \sum _{\gamma
\in \Gamma} \frac{\frac{\pi}{2}\xi _r(\gamma ,z)} {e^{2\nu
(z)}\Delta \vp(z)+\tau_{\psi}(z)},
\end{eqnarray}
and
\begin{eqnarray}\label{d-}
D^-_f(\Gamma) := \liminf _{r \to R_X} \inf _{z\in X} \sum _{\gamma
\in \Gamma} \frac{\frac{\pi}{2}\xi _r(\gamma ,z)} {e^{2\nu(z)}\Delta
\vp(z)},
\end{eqnarray}
where $\xi _r$ is defined as in section \ref{intro}

\subsection{A sub-mean value lemma}

We will have occasion to use the following result.

\begin{lem}\label{sub-mvt}
For any function $F$,
\begin{eqnarray*}
\int _X F(w) \xi _r (z,w) e^{-2\nu (w)}\ii dw \wedge d\bar w &=&
\frac{1}{c_r} \int _{D_r(z)} F f(\rho _z) d\rho _z \wedge *d\rho _z \\
&=& \frac{1}{c_r}\int _0 ^r t f(t)\left (\int _{S_t(z)}F * dE_z
\right )dt.
\end{eqnarray*}
Thus, in view of (\ref{sub-pois}) of Lemma \ref{mvt}, if $h$ is
subharmonic then
\begin{eqnarray}
h(z) &\le& \int _X \xi _r (z,w)h(w) e^{-2\nu (w)} \ii dw \wedge
d\bar w \label{area-mvp}
\end{eqnarray}
with equality if $h$ is harmonic.
\end{lem}

\subsection{Interpolation and sampling theorems}

\begin{thm}\label{interp-thm}
Let $X$ be a finite open Riemann surface with metric
$g=e^{-2\psi}|dz|^2$ and let $\vp$ be a weight function on $X$ such
that, for some $c>1$, $\frac{1}{c} \le e^{2\nu} \Delta \vp \le c$,
\begin{eqnarray}\label{fund-hyp-h}
e^{2\nu}\Delta (\vp) +\tau _{\psi} \ge \frac{1}{c} \quad and \quad
e^{2\nu}|\di u_{\psi}|^2 \le c.
\end{eqnarray}
Then every uniformly separated sequence $\Gamma$ satisfying
$D^+_f(\Gamma) < 1$ is an interpolation sequence.
\end{thm}

\begin{thm}\label{sampling-thm}
Let $X$ be a finite open Riemann surface with metric $g=e^{-2\psi}
|dz|^2$ and $\Gamma \subset X$ a uniformly separated sequence.
Suppose $\vp$ is a weight function on $X$ such that, for some $C>1$,
$\frac{1}{C} \le e^{2\nu} \Delta \vp \le C$.  Assume also that the
metric $g$ is bounded above by the fundamental metric $e^{-2\nu}$
(i.e., $u_{\psi} \ge 0$) and moreover satisfies the differential
inequality
\begin{eqnarray}\label{diff-ineq}
\quad 2e^{2\nu} |\di u_{\psi}|^2 \le e^{2\nu} \Delta u_{\psi}.
\end{eqnarray}
Then $\Gamma$ is a sampling sequence whenever $D^-_f(\Gamma) > 1$.
\end{thm}

Theorems \ref{interp-thm} and \ref{sampling-thm} imply Theorem
\ref{main-thm}.

\subsection{Existence of metrics satisfying (\ref{diff-ineq})}

Let $X$ be an open Riemann surface. Observe that a function $u$ on
$X$ satisfies (\ref{diff-ineq}) if and only if
$$\Delta  (-e^{-2u}) \ge 0.$$  Since the function $-e^{-2u}$ is
bounded above, we immediately obtain the following proposition.

\begin{prop}\label{no-u}
Let $X$ be a compact or parabolic Riemann surface.  Then any
function satisfying the inequality {\rm (\ref{diff-ineq})} is
constant.
\end{prop}

\noi In particular, we may assume that when $X$ is parabolic, the
metric $g$ in Theorem \ref{sampling-thm} is just the fundamental
metric. \vskip0.25in

Let us turn now to the hyperbolic case and suppose that $o \in X$.
Then
$$
\rho _o \Delta \rho _o  = \rho _o ^2 \Delta E(o,\cdot ) + \rho _o ^2
|\di E(o,\cdot)|^2  = \rho _o ^2 |\di E(o,\cdot)|^2 = |\di \rho
_o|^2,
$$
and thus
\begin{eqnarray*}
\Delta \rho _o ^2 = 4|\di \rho _o|^2 = |d \rho _o |^2.
\end{eqnarray*}
We let
$$
u = -\frac{1}{2} \log (1 - \rho _o ^2).
$$
The reader versed in Several Complex Variables will recognize this
function as the {\it negative log-distance-to-the-boundary}.
Calculating, we have
\begin{eqnarray*}
\Delta u - 2 |\di u|^2 &=& \left (\frac{\Delta \rho _o^2 }{2(1-\rho
_o ^2)}+ \frac{|\di \rho _0 ^2|^2}{2(1-\rho _o ^2)^2}\right ) - 2
\frac{|\di \rho _o
^2|^2} {4(1-\rho _0^2)^2}\\
&=& \frac{|d \rho _0 |^2}{2(1-\rho _0 ^2)}\ge 0.
\end{eqnarray*}
Moreover, observe that
$$
|\di u|^2 = \frac{\rho _o^2 |\di \rho _o|^2}{(1-\rho _o ^2)^2}.
$$
\begin{prop}\label{some-u}
Let $X$ be a hyperbolic Riemann surface.  Then there is always a
non-trivial function $u \ge 0$ satisfying the differential
inequality {\rm (\ref{diff-ineq})}. Moreover, if $X$ is a finite
hyperbolic surface, then one can choose $u$ such that $e^{2\nu}|\di
u|^2$ is uniformly bounded.
\end{prop}

\noi {\it Sketch of proof.} It remains only to verify the last
assertion.  By compactness, it suffices to prove the desired
estimate in a neighborhood of the form $\{ z \in \C \ ;\ |z|<1,\ \im
z >0\}$ in the upper half plane.  As above, one can take $\nu = -
\log {\im z}$ the Poincar\' e potential.  Moreover, one can show
that
$$1-\rho_o(z) = \im z + \ {\rm higher\ order\ terms}.$$
The proposition now follows.\qed

\section{A Theorem of Ohsawa on the solution of $\dbar$}\label{dbar}

In our proof of the interpolation theorem, we require a theorem for
solving $\dbar$ with certain $L^2$ estimates.  Such a theorem has
been stated by Ohsawa \cite{o} in a very general situation, but
there seem to be counterexamples at this level of generality (see
\cite{siu2}). However, in the case of Riemann surfaces there is a
short proof of Ohsawa's theorem.  Since it is not easily accessible
in the literature, we shall give a proof here using methods adapted
from \cite{siu2}.

Let $X$ be a Riemann surface with conformal metric $g=e^{-2\psi}$
and let $V\to X$ be a holomorphic line bundle with Hermitian metric
$h=e^{-2\xi}$ that is allowed to be singular, i.e., $\xi$ may be in
$L^1 _{\ell oc}$.  One has the Bochner-Kodaira identity
\begin{eqnarray}\label{bk-basic}
||\dbar ^* \!\! \beta||^2 = ||\nbar \beta ||^2 + (2 e^{2\psi} \Delta
(\xi + \psi)  \beta, \beta),
\end{eqnarray}
where
$$
\nbar (f d\bar z ):=(f_{\bar z}+2\psi _{\bar z}f)d\bar z^{\otimes
2}.
$$
Indeed, straight-forward calculations show that the formal adjoints
$\dbar ^*$ of $\dbar$ and $\nbar ^*$ of $\nbar$ are given by
\begin{eqnarray}\label{formal-adjoint}
\dbar ^*(h d\bar z) = -e^{2\psi}\left ( \frac{\di h}{\di z}-2
\frac{\di \xi}{\di z} h \right )\ {\rm and}\  \nbar ^* (hd\bar
z^{\otimes 2}) = - e^{2\psi} \left ( \frac{\di h}{\di z} -
2\frac{\di \xi} {\di z} h \right ).
\end{eqnarray}
Using these, another calculation shows that $\dbar \dbar ^*\beta -
\nbar ^* \nbar \beta = 2 e^{2\psi} \Delta (\xi + \psi) \beta$, which
gives (\ref{bk-basic}).

We shall now make a simple but far-reaching modification of the
identity (\ref{bk-basic}).  To this end, let $e^{-2\xi}$ and
$e^{-2(\xi -u)}$ be two metrics of the same line bundle. (Thus $u$
is a globally defined function.)  We assume moreover that $e^{2(\psi
-u)}|\di u |^2$ is uniformly bounded.

Formula (\ref{formal-adjoint}) implies that
\begin{eqnarray}\label{compare-f-as}
\dbar ^* _{\xi - u} \beta = \dbar ^* _{\xi} \beta - 2e^{2\psi}\di u
\wedge \beta .
\end{eqnarray}
Substituting (\ref{formal-adjoint}) into the Bochner-Kodaira
identity (\ref{bk-basic}), we obtain
\begin{eqnarray}\label{bk-2w}
||\dbar ^* _{\xi - u} \beta||^2_{\xi} &=& ||\nbar \beta ||^2_{\xi}  \\
\nonumber && \qquad + \left ( 2 e^{2\psi} \left \{ \Delta (\xi +
\psi) - 2 |\di u|^2
\right \}  \beta, \beta \right )_{\xi}  \\
\nonumber && \qquad \qquad \qquad \qquad - 2 {\rm Re} (\dbar ^*
_{\xi} \beta , 2e^{2\psi} \di u \wedge \beta ) _{\xi}.
\end{eqnarray}
Identity (\ref{bk-2w}) is sometimes called the Bochner-Kodaira
identity with two weights. The Cauchy-Schwarz inequality then shows
that for any $\ve > 0$ we have
\begin{eqnarray}\label{pre-solution-est}
(1+\ve ^{-1}) ||\dbar ^* _{\xi -u} \beta||^2_{\xi} \ge \left ( 2
e^{2\psi} \left \{ \Delta (\xi + \psi) - 2 (1+ \ve) |\di u|^2 \right
\}  \beta, \beta \right )_{\xi}.
\end{eqnarray}
Letting $Tf:=\dbar (e^{-u} f)$, we can rewrite
(\ref{pre-solution-est}) as
\begin{eqnarray}\label{solution-est}
||T^* \beta ||^2 _{\xi -u} \ge C _{\ve} \left ( 2 e^{2(\psi -u)}
\left \{ \Delta (\xi + \psi) - 2 (1+ \ve) |\di u|^2 \right \} \beta,
\beta \right )_{\xi -u}.
\end{eqnarray}
Suppose now that, for some $\delta >0$, one has the estimate
$$
e^{2(\psi -u)}\left ( \Delta (\xi + \psi) - 2 |\di u|^2 \right ) \ge
\delta.
$$
Since $e^{2(\psi -u)}|\di u|^2$ is bounded, we may choose $\ve >0$
sufficiently small in (\ref{solution-est}) to obtain
\begin{eqnarray}\label{next-one}
||T^* \beta ||^2 _{\xi -u} \ge C ||\beta||_{\xi -u}.
\end{eqnarray}
A standard Hilbert space argument yields a function $f$ such that
$$Tf = \alpha$$
with the estimate
\begin{eqnarray}\label{l2-est-for-t}
\int _X |f|^2 e^{2u} e^{-2\xi} dA_g \le C ||\alpha||_{\xi -u}.
\end{eqnarray}
Finally, choosing $u = u_{\psi} = \psi - \nu$, $\vp := \xi - 2u$ and
$U= e^{-u} f$ gives the following.

\begin{thm}\label{o-dbar-thm}{\rm \cite{o}}
Suppose that for some $\delta > 0$,
$$
e^{2\nu} \Delta \vp + \tau _{\psi} \ge \delta \quad {\rm and} \quad
e^{2\nu}|\di u_{\psi}|^2 < \frac{1}{\delta}.$$ Then there exists a
constant $C=C_{\delta}$ such that for any $\alpha$ satisfying
$$
\int _X \ii \alpha \wedge \bar \alpha e^{-2u_{\psi}} e^{-2\vp} <
+\infty,
$$
the equation $\dbar U = \alpha$ has a solution satisfying
\begin{eqnarray*}
\int _X |U|^2 e^{-2\vp} dA_g \le C \int _X \alpha \wedge \bar \alpha
e^{-2u_{\psi}} e^{-2\vp}.
\end{eqnarray*}
\end{thm}

\section{Compact Riemann surfaces}\label{cpct}

At this point we take a short detour to consider the problems of
sampling and interpolation on elliptic Riemann surfaces. While the
essence of this situation is different from that of open Riemann
surfaces, we note that the estimates on the solution of the
$\overline{\partial}$ problem discussed in the previous section are
applicable.

Let $X$ then be a compact Riemann surface and let $V \to X$ be a
holomorphic line bundle.  We denote by $V_x$ the fiber of $V$ over
$x \in X$. Then $\Gamma$ is interpolating if and only if the
evaluation map
\begin{eqnarray}\label{eval-map}
H^0(X,L) \ni s \mapsto  \sum _{\gamma \in \Gamma} s(\gamma) \in
\bigoplus _{\gamma \in \Gamma} V_{\gamma}
\end{eqnarray}
is surjective, and sampling if and only if (\ref{eval-map}) is
injective.

Let $\Lambda$ be the line bundle corresponding to the effective
divisor $\Gamma$.  One can understand the situation completely using
the short exact sequence of sheaves
\begin{eqnarray*}
0 \to \sco _X(L\tensor \Lambda ^*) \to \sco _X(L) \to \bigoplus
_{\gamma \in \Gamma} \cv _{\gamma} \to 0,
\end{eqnarray*}
where $\cv _{\gamma} (U) = V_{\gamma}$ if $\gamma \in U$ and $\cv
_{\gamma} (U) =0$ if $\gamma \not \in U$.  Passing to the long exact
sequence, we have that
$$
0 \to H^0(X,L\tensor \Lambda ^*) {\buildrel i_0 \over
\longrightarrow} H^0(X,L) {\buildrel e_{\Gamma} \over
\longrightarrow} \bigoplus _{\gamma \in \Gamma} V_{\gamma}
{\buildrel \delta _0 \over \longrightarrow} H^1 (X,L\tensor \Lambda
^*) {\buildrel i_1 \over \longrightarrow} H^1 (X,L) \to ...
$$
We see that $e$ is injective if and only if Image$(i_0) = \{0\}$ and
surjective if and only if $i_1$ is injective, i.e., Image$(\delta
_0)=\{0\}$.  We then have the following proposition.
\begin{prop} \label{alg}
Let $X$ be a compact Riemann surface of genus $g$, $\Gamma \subset
X$ a finite subset and $L \to X$ a holomorphic line bundle.
\begin{enumerate}
\item If $\# \Gamma < {\rm deg}(L) +2 -2g$, then $\Gamma$ is
interpolating. \item If $\# \Gamma > {\rm deg}(L) $, then $\Gamma$
is sampling.
\end{enumerate}
\end{prop}

\pf To establish 1, note that by Serre duality, $h^1(X,L\tensor
\Lambda ^*)  = h^0(X,K_X \tensor \Lambda \tensor L^*)$, and the
latter vanishes if
$$
\# \Gamma + 2g-2 - {\rm deg}(L) = {\rm deg}(K_X \tensor \Lambda
\tensor L^*) < 0.
$$
Similarly, if ${\rm deg}(L) - \# \Gamma = {\rm deg}(L \tensor
\Lambda ^*) <0$, then $h^0(X,L\tensor \Lambda ^*) =0$.\qed

\medskip

\noi Part 1 of Proposition \ref{alg} can also be proved using
Theorem \ref{o-dbar-thm}.  Because it is similar to the proof of our
main interpolation theorem, we sketch this method here.

\medskip

\noi {\it Analytic proof of Proposition \ref{alg}.1.} Let $\sum
v_{\gamma} \in \bigoplus V_{\gamma}$. First, observe that there is a
smooth section $\eta$ of $L$ such that $\eta (\gamma) = v_{\gamma}$
for all $\gamma \in \Gamma$. In fact, by the usual cutoff method, we
can take $\eta$ supported near $\Gamma$ and holomorphic in a
neighborhood of $\Gamma$.

Fix a conformal metric $e^{-2\psi}$ on $X$. Let $\tau$ be the
canonical section of $\Lambda$ corresponding to the divisor
$\Gamma$. By the degree hypothesis, there is a metric $e^{-2\vp}$
for the line bundle $L\tensor \Lambda ^*$ such that the curvature
$\ii \di \dbar (\vp + \psi)$ of $L\tensor \Lambda ^* \tensor K_X^*$
is strictly positive on $X$. Then $e^{-2\vp}/|\tau|^2$ is a singular
metric for $L$ such that the curvature current of $e^{-2(\vp
+\psi)}/|\tau|^2$ is still strictly positive on $X$. Moreover, since
$\eta$ is holomorphic in a neighborhood of $\Gamma$, we have $\int
_X |\dbar \eta |^2 |\tau|^{-2}e^{-2\vp} < + \infty.$ By Theorem
\ref{o-dbar-thm} (with $u_{\psi} \equiv 0$; {\it c.f.} Proposition
\ref{no-u}) there is a section $u$ of $L$ such that $\dbar u = \dbar
\eta$ and $\int _X |u|^2 |\tau|^{-2}e^{-2(\vp+\psi)}< + \infty.$ But
since $\tau$ vanishes on $\Gamma$, so does $u$.  Thus $\sigma = \eta
- u$ is holomorphic and solves the interpolation problem. \qed

\medskip

\rmk {\it We note that if $e^{-2\vp}$ is a metric for a holomorphic
line bundle $L$, then we have
$${\rm deg}(L) = \frac{1}{2\pi}\int _X \Delta \vp,$$
showing the resemblance between Proposition \ref{alg} and our main
theorems.}

\section{Proofs of Theorems \ref{interp-thm} and \ref{sampling-thm}}
\label{proofs}

\subsection{Functions and singular weights}

In this paragraph we define certain functions that play important
roles in the proofs of Theorems \ref{interp-thm} and
\ref{sampling-thm}.

\subsubsection*{A local construction of a holomorphic function}

In the proofs of Theorems \ref{interp-thm} and \ref{sampling-thm}
we will need, for each $\gamma \in \Gamma$, a holomorphic function
defined in a neighborhood of $\gamma$ and satisfying certain
global estimates. For reasons that will become clear later, the
size of this neighborhood cannot be taken too small. As a
consequence, we must overcome certain difficulties presented by
the topology of the neighborhood.

\begin{lem}\label{local-fn}
Let $X$ be a finite open Riemann surface. Assume $e^{2\nu}\Delta
\vp$ is bounded above and below by positive constants. Let $\Gamma$
be a uniformly separated sequence and $\sigma = \sigma (\Gamma)$.
There exists a constant $C=C_{\Gamma}>0$ and, for each $\gamma \in
\Gamma$, a holomorphic function $F_{\gamma} \in \sco
(D_{\sigma}(\gamma))$ such that $F_{\gamma} (\gamma) = 0$ and for
all $z \in D_{\sigma} (\gamma)$,
\begin{eqnarray}\label{local-est}
\frac{1}{C} e^{-2\vp (\gamma)} \le  \left | e^{-2\vp +
2F_{\gamma}} \right | \le C e^{-2\vp (\gamma)}.
\end{eqnarray}
\end{lem}

\pf Let $G$ be the Green's function for the domain
$D_{2\sigma}(\gamma)$.  Consider the function
$$
T_{\gamma} (z):= \frac{2}{\pi} \int _{D_{2\sigma}(\gamma)}
\!\!\!\!\!\!\!\! -G(z,\zeta) \Delta \vp (\zeta).
$$
By Green's formula, we have that
$$
T_{\gamma} (z) = - \vp (z) + \frac{1}{2\pi} \int
_{S_{2\sigma}(\gamma)} \!\!\!\!\!\!\!\! \vp (\zeta) * d_{\zeta}G
(z,\zeta).
$$
We claim that the harmonic function
$$
h_{\gamma} := \frac{1}{2\pi}\int _{S_{2\sigma}(\gamma)}
\!\!\!\!\!\!\!\! \vp (\zeta)
* d_{\zeta}G (z,\zeta)
$$
has a harmonic conjugate, i.e., it is the real part of a
holomorphic function.  Indeed, if $\cc$ is a Jordan curve in
$D_r(\gamma)$, then
\begin{eqnarray}\label{cohomology}
\int _{\cc} *dh_{\gamma} (z) &=& \frac{1}{2\pi} \int _{S_{2\sigma}
(\gamma)} \!\!\!\!\!\!\!\! \vp (\zeta) * d_{\zeta} \left ( \int
_{\cc} *d_z G (z,\zeta) \right ).
\end{eqnarray}
Since $S_{2\sigma}(\gamma) \cap \cc = \emptyset $, the
function $z \mapsto G(z,\zeta)$ is harmonic and thus $*d_z G
(z,\zeta)$ is a closed form.  It follows that the term in the
parentheses on the right hand side of (\ref{cohomology}) depends
only on the homology class $[\cc]\in H_1(X,\Z)$. Since $H_1(X,\Z)$
is discrete and $*d_z G (z,\zeta)$ is continuous in $\zeta$, we
see that the right hand side of (\ref{cohomology}) vanishes, as
claimed.

Let
$$
H_{\gamma} := h_{\gamma} + \ii \int _{\gamma} ^z \!\!\!\! *dh
_{\gamma}
$$
be the holomorphic function whose real part is $h_{\gamma}$, and
let $F_{\gamma} := H_{\gamma} - H_{\gamma} (\gamma).$  We have
\begin{eqnarray*}
\left | 2 \vp (\gamma ) -2\vp (z) + 2\re F_{\gamma} (z) \right | =
2 \left | T_{\gamma} (\gamma) - T_{\gamma} (z) \right | \le
2|T_{\gamma} (\gamma)|+ 2|T_{\gamma} (z)|.
\end{eqnarray*}
Taking exponentials and applying Lemma \ref{hh} completes the
proof.\qed

\subsubsection*{A function with poles along $\Gamma$}

For $z,\zeta\in X$ and $r<R_X$, let
\begin{eqnarray*}
I(\zeta,z) &=& \int_X \xi_r(\zeta,w)E(w,z)e^{-2\nu(w)}\ii dw \wedge d\bar w \\
&=& \frac{1} {c_r}\int _a ^r tf(t) \left ( \int_
{S_t(\zeta)}E(w,z) *dE_{\zeta} (w) \right ) dt.
\end{eqnarray*}
Since $E$ is a fundamental solution to the Laplacian,
$$
e^{2\nu(z)}\Delta_zI(\zeta,z)=\frac{\pi}{2}\int_X\xi_r(\zeta,w)
\delta_z(w)=\frac{\pi}{2} \xi_r(\zeta,z).
$$
Next it follows from (\ref{area-mvp}) that, since $E(\cdot, z)$ is
subharmonic, $E(\zeta,z)\leq I(\zeta,z)$ and, since $E(\cdot, z)$
is harmonic in the region $\{w\in X:\rho_\zeta(w)>r\}$,
$E(\zeta,z)=I(\zeta,z)$ if $\rho_z(\zeta)>r$. Moreover, in view of
(\ref{h=e}), an application of (\ref{pois}) shows that
$$
\frac{1}{2\pi} \int_{S_t(\zeta)}\!\!\!\!\!\!\!\!
E(w,z)*dE_{\zeta}(w)=E(z,\zeta)-{\bf 1}_{D_t(z)}(\zeta)
\left(E(z,\zeta)-\log t\right).
$$
We see that
\begin{eqnarray*}
I(\zeta,z)&=& \frac{2\pi}{c_r} \left ( \log (\rho _z(\zeta)) \int
_0 ^{\rho _z(\zeta)} \!\!\!\!\!\!\!\!\!\!\! tf(t) dt + \int _{\rho
_z(\zeta)} ^r \!\!\!\!\!\!\!\!\! tf(t)\log t dt\right )
\end{eqnarray*}
if $\rho_z(\zeta)<r$. Note that
$$
\left | \frac{1}{c_r}\int_{\rho _z (\zeta)} ^r \!\!\!\! tf(t) \log
(t) dt \right | \le D_r,
$$
where $D_r$ depends only on $r$.  We then have
\begin{eqnarray*}
|I(\zeta,z)|\leq K_r
\rho_z(\zeta)\left|\log(\rho_z(\zeta))\right|+ D_r
\end{eqnarray*}
for all $z,\zeta\in X$ satisfying $\rho_z(\zeta)<r$. Since the
expression on the right hand side is bounded by a constant that
depends only on $r$, we have
\begin{eqnarray}\label{i-est}
|I(\zeta,z)|\leq C_r
\end{eqnarray}
whenever $\rho_z(\zeta)<r$.

Let $\Gamma$ be a discrete sequence.  We define the function
$$v_r(z)=\sum_{\gamma\in\Gamma}\left(E(\gamma,z)-I(\gamma,z)\right).$$
By the preceding remarks, $v_r(z)\leq0$ and
$$
v_r(z)=\sum_{\gamma\in\Gamma\cap D_r(z)} \left ( E(\gamma,z) -
I(\gamma ,z) \right).
$$
Moreover,
\begin{eqnarray}\label{lap-v}
e^{2\nu} \Delta v_r = \frac{\pi}{2}\sum _{\gamma\in\Gamma}
(e^{2\nu} \delta _{\gamma } - \xi _r (\gamma, \cdot)).
\end{eqnarray}
Writing
$$
X_{\Gamma ,\ve}  := \left \{ z\in X\ ;\ \min _{\gamma \in \Gamma}
\rho _{\gamma} (z) > \ve \right \},
$$
we have the following lemma.
\begin{lem}\label{v-estimate}
Let $\Gamma$ be a sparse, uniformly separated sequence and let $\ve
\le \sigma (\Gamma)$.  The function $v_r$ is uniformly bounded on
$X_{\Gamma ,\ve}$. Moreover, $v_r$ satisfies the following estimate:
if $\gamma \in \Gamma$ and $\rho _{\gamma } (z) < \sigma(\Gamma)$,
then
\begin{eqnarray}\label{v-pole}
\left | v_r(z) - \log \rho _{\gamma} (z) \right | \le C_{r,\ve}.
\end{eqnarray}
\end{lem}

\pf Let $z\in X_{\Gamma,\ve}$. Since $\Gamma$ is sparse, there are
at most $N=N_{r,0}$ members of $\Gamma$, say $\gamma_{1},\dots,
\gamma_{N}$, lying in $D_r(z)$, and so
$$
|v_r(z)|\leq \sum_{j=1}^{N} (|E(\gamma_{j},z)| + |I(\gamma
_{j},z)|) \leq \sum_{j=1}^{N}\left ( \left| \log (\rho_z (\gamma_
{j}) )\right|+C_r\right).
$$
Note that the number $N$ does not depend on $z$.  Since $\ve <
\rho_z(\gamma_j) < r$, the term involving the logarithm has a bound
that depends only on $\ve$ and $r$. We thus see that $v_r$ is
uniformly bounded on $X_{\Gamma,\ve}$.

Let $\gamma\in\Gamma$.  Since $\Gamma$ is sparse, there are at most
$N=N_{r,\ve}$ elements of $\Gamma$ that lie in $D_r(D_\ve(\gamma))$.
We write $\Gamma\cap
D_r(D_\ve(\gamma))=\{\gamma_{1},\dots,\gamma_{N} \}$, where
$\gamma_{1}=\gamma$. Again, $N$ does not depend on $z$. Then
$$
|v_r(z)-\log\rho_z(\gamma)| \leq \left ( \sum_{j=2}^N
|E(\gamma_{j},z)| + \sum_{j=1}^N |I(\gamma_{j},z)| \right ) +
|E(\gamma,z) - \log\rho_z (\gamma)|.
$$
The first sum is bounded because $\sigma(\Gamma) < \rho_z
(\gamma_{j}) < r$ for $j=2,\dots,N$. The second sum is bounded by
(\ref{i-est}), and the third term vanishes.  This completes the
proof of the lemma.\qed

\subsubsection*{A function with bumps along $\Gamma$}

In this paragraph, we shall use area forms associated to the
points of $\Gamma$.  We define
$$dA_{E,\gamma}(\zeta) := d \rho _{\gamma}(\zeta) \wedge *d\rho _{\gamma}
(\zeta).$$

Let
$$
A_{E,\gamma} (D) := \int _D dA_{E,\gamma}.
$$
Given a distribution $f$, we consider its regularization
$$
\frac{1}{A_{E,\gamma}(D_{\ve}(z))} \int _{D_{\ve}(z)} f
dA_{E,\gamma}
$$
using the area element $dA_{E,\gamma}$, where $\gamma \in \Gamma$.

Observe that
\begin{eqnarray} \nonumber
A_{E,\gamma} (D_{\ve}(\gamma)) &=& \int _{D_{\ve}(\gamma)} d\rho
_{\gamma} \wedge *d\rho _{\gamma} = \int _{D_{\ve}(\gamma)} \rho
_{\gamma} d\rho _{\gamma} \wedge *dE(\gamma,\cdot) \\
\label{exact-area} &=& \int _0 ^{\ve} t \left (\int _{S_t(\gamma)}
*dE(\gamma , \cdot) \right ) dt = 2\pi \int _0 ^{\ve} t\ dt = \pi
\ve ^2.
\end{eqnarray}

Consider the function
\begin{eqnarray*}
v _{r,\ve}(z) &=& t \sum _{\gamma \in \Gamma} \frac{1} {\pi \ve ^2}
\int _{D_{\ve}(\gamma)} \left ( E(\zeta,z) - I(\zeta,z) \right )
dA_{E,\gamma}(\zeta)
\end{eqnarray*}
where $0<<t<1$.

\begin{lem}\label{v-estimates} The function
$v_{r,\ve}$ has the following properties.
\begin{enumerate}
\item
\begin{eqnarray*}
e^{2\nu (z)} \Delta v_{r,\ve}(z) &=& t \sum _{\gamma \in \Gamma}
 \frac{1}{2\ve^2}e^{2\nu (z)}|d \rho _{\gamma}(z)|^2  {\bf
1}_{D_{\ve}(z)} \\
&& \qquad - t \sum _{\gamma \in \Gamma} \frac{1}{2\ve^2} \int
_{D_{\ve}(\gamma)}\xi _r (\cdot ,z) dA_{E,\gamma}.
\end{eqnarray*}
In particular,
$$
\lim _{\ve \to 0} e^{2\nu} \Delta v_{r,\ve} =\frac{\pi}{2} t \sum
_{\gamma \in \Gamma} \left ( e^{2\nu} \delta _\gamma-\xi _r
(\gamma,\cdot) \right )$$ in the sense of distributions. \item There
exists a positive constant $C_{r,\ve}$ such that
\begin{eqnarray} \label{1st-v-estimate}
z \in X &\Rightarrow& -C_{r,\ve} \le v _{r,\ve} (z) \le 0
\end{eqnarray}
and for any $\gamma \in \Gamma$,
\begin{eqnarray}\label{2nd-v-estimate}
\rho _\gamma (z) < \ve \!\!\!\! & \Rightarrow & \!\!\!\! \left | \ v
_{r,\ve} (z) - \frac{t}{\pi \ve ^2} \int _{D_{\ve} (\gamma) }
\!\!\!\!\!\!\!\! E(\zeta , z) dA_{E,\gamma}(\zeta) \right | \le
C_{r,\ve}
\end{eqnarray}
\end{enumerate}
\end{lem}

\pf {\it 1.} The formula for the Laplacian is a straightforward
calculation, and the limit is a standard consequence of the
regularization of currents.

\medskip

\noi {\it 2.} Since $E(\zeta,z)=I(\zeta,z)$ whenever $\rho
_z(\zeta) > r$, we have, in view of formula (\ref{exact-area}),
$$
v_{r,\ve}(z) = \sum _{\gamma \in D_{\ve}(D_r(z))} \frac{t} {\pi \ve
^2} \int _{D_{\ve}(\gamma)} \!\!\!\!\!\!\!\! \left ( E(\zeta,z) -
I(\zeta,z) \right ) dA_{E,\gamma}(\zeta).
$$
Choose $\gamma \in \Gamma$. Since $\Gamma$ is sparse, there exist
$\gamma _1,...,\gamma _N \in \Gamma - \{\gamma\}$ such that for all
$z \in D_{\ve} (\gamma)$
\begin{eqnarray*}
&& v_{r,\ve}(z) = \frac{t} {\pi \ve ^2} \int _{D_{\ve}(\gamma)}
\!\!\!\!\!\!\!\! \left (E(\zeta, z) -
I(\zeta,z) \right ) dA_{E,\gamma} \\
&& \qquad  \qquad  \qquad  + \sum _{j=1} ^N \frac{t} {\pi \ve ^2}
\int _{D_{\ve}(\gamma _j)} \!\!\!\!\!\!\!\! \left (E(\zeta, z) -
I(\zeta,z) \right ) dA_{E,\gamma _j}
\end{eqnarray*}
Moreover, $N$ is independent of $\gamma$, and depends only on $r$
and $\ve$.  It follows that
\begin{eqnarray*}
&& \left | v_{r,\ve} (z) - \frac{t}{\pi \ve ^2}\int _{D_{\ve}
(\gamma)} \!\!\!\!\!\!\!\! E(\cdot ,z)
dA_{E,\gamma} \right | \\
&& \qquad \le  \frac{t} {\pi \ve ^2} \int _{D_{\ve}(\gamma)}
\!\!\!\!\!\!\!\! |I(\cdot , z)| dA_{E,\gamma} + \sum _{j=1} ^N
\frac{t}{\pi \ve ^2} \int _{D_{\ve}(\gamma _j)} \left ( |E(\cdot ,
z)| + |I(\cdot , z)| \right ) dA_{E,\gamma _j}.
\end{eqnarray*}
We have estimates for $I(\zeta,z)$ as in the proof of Lemma
\ref{v-estimate}, and since, by uniform separation, $\rho _z(\zeta)
> \sigma$ for any $\zeta \in D_{\ve} (\gamma _j)$, we can estimate
the right hand side by a constant that depends only on $r$.   This
proves (\ref{2nd-v-estimate}), and (\ref{1st-v-estimate}) follows
from (\ref{2nd-v-estimate}), Lemma \ref{v-estimate} and the fact
that $v_r \le 0$.\qed

\medskip

Finally, we shall have use for the following lemma.

\begin{lem}\label{min}
For any $z \in D_{\ve} (\gamma)$,
\begin{eqnarray}\label{e-ineq}
\frac{1}{A_{E,\gamma} (D_{\ve}(\gamma))} \int _{D_{\ve}(\gamma)}
\!\!\!\!\!\!\!\! E(z,\zeta) dA_{E,\gamma}(\zeta) \le \log
\frac{1}{\ve} +\frac{1}{2}.
\end{eqnarray}
\end{lem}

\pf Observe that if $z \in D_{\ve}(\gamma)$ and $t\in (0,\ve]$, then
$$
\int _{S_t(\gamma)} *d_{\zeta} E(z,\zeta) = \int _{D_t(\gamma)}
d_{\zeta}*d_{\zeta} E(z,\zeta) = 2\pi {\bf 1} _{D_t(\gamma)}(z)
\le 2\pi.
$$
Applying Green's formula (\ref{grn-fa}) with $f=E(z,\cdot)$ and
$h= E(\gamma,\cdot)$, we obtain
$$
\int _{S_t(\gamma)} E(z,\zeta) *d_{\zeta} E(\gamma,\zeta) = \int
_{S_t(\gamma)} E(\gamma,\zeta) *d_{\zeta} E(z,\zeta).
$$
We thus have
\begin{eqnarray*}
- \int _{D_{\ve}(\gamma)}\!\!\!\!\!\!\!\! \log \rho _z d\rho
_{\gamma} \wedge *d\rho _{\gamma} &=& -\int _0 ^{\ve} t \left ( \int
_{S_t(\gamma)} E(z,\zeta)
*d_{\zeta} E(\gamma,\zeta) \right ) dt \\
&=& -\int _0 ^{\ve} t \left ( \int _{S_t(\gamma)} E(\gamma,\zeta)
*d_{\zeta} E(z,\zeta) \right ) dt \\
&=& -\int _0 ^{\ve} t \log t \left ( \int _{S_t(\gamma)}
*d_{\zeta} E(z,\zeta) \right ) dt \\
&\le & -2\pi \int _0 ^{\ve} t \log t dt = \pi \ve ^2 \left (
\frac{1}{2} - \log \ve \right ).
\end{eqnarray*}
The lemma now follows from (\ref{exact-area}).\qed

\subsection{Proof of Theorem
\ref{interp-thm}}\label{interp-section}

Let $(s_{\gamma}) \in \fb ^2 (\Gamma , g,\vp)$.  We begin by
constructing a smooth function $\eta \in L^2(X,g,\vp)$ that
interpolates $(s_{\gamma})$. To this end, let $\chi \in \cc
^{\infty} _0 ([0,\sigma))$ satisfy
$$
0 \le \chi \le 1,\quad \chi | [0,\sigma /2] \equiv 1 \quad {\rm
and} \quad |\chi '| \le \frac{3}{\sigma}.
$$
We define
$$
\eta (z) := \sum _{\gamma \in \Gamma} \chi \circ \rho _{\gamma}
(z) s_{\gamma} e^{F_{\gamma}(z)},
$$
where $F_{\gamma}$ is as in Lemma \ref{local-fn}.  Observe that
$\eta (\gamma ) = s_{\gamma}$ for all $\gamma \in \Gamma$, and
that
\begin{eqnarray*}
\int _X |\eta|^2 e^{-2\vp } dA_g &=&  \sum _{\gamma \in \Gamma}
|s_{\gamma}|^2 \int _{D_{\sigma} (\gamma)} | \chi \circ \rho
_{\gamma} |^2 \left | e^{2F_{\gamma} -
2\vp} \right | dA_g\\
&\le & C \sum _{\gamma \in \Gamma} |s_{\gamma}|^2 e^{-2\vp
(\gamma)} A_g(D_{\sigma}(\gamma)) < +\infty.
\end{eqnarray*}

Next we wish to correct $\eta$ by adding to it a function $U$ that
lies in $L^2 (X,g,\vp)$ and vanishes along $\Gamma$.  The standard
approach is to solve the equation $\dbar U = \dbar \eta$ with
singular weights, using Ohsawa's $\dbar$ Theorem \ref{o-dbar-thm}.
The singular weight we will use is the weight
$$\tilde \vp := \vp + v_r,$$
and one computes that
\begin{eqnarray}\label{dbar-eta}
\dbar \eta = \sum _{\gamma \in \Gamma} \chi ' (\rho _{\gamma})
\dbar \rho _{\gamma} s_{\gamma} e^{F_{\gamma}}.
\end{eqnarray}
Since $D^+_f(\Gamma) < 1$, there exist $r < R_X$ and $\delta
> 0$ such that
\begin{eqnarray*}
e^{2\nu} \Delta \tilde \vp + \tau _{\psi} &=& e^{2\nu} \Delta
\vp + \tau _{\psi} + e^{2\nu} \Delta v_r \\
&\ge& \left ( e^{2\nu} \Delta \vp + \tau _{\psi}\right ) \left ( 1
- \sum _{\gamma \in \Gamma} \frac{\frac{\pi}{2}\xi _r (\cdot
,\gamma)}{e^{2\nu} \Delta \vp + \tau _{\psi}} \right ) \\
&>& \delta (e^{2\nu} \Delta \vp + \tau _{\psi}),
\end{eqnarray*}
where the first inequality follows from (\ref{lap-v}).  It follows
from hypothesis (\ref{fund-hyp-h}) in Theorem \ref{interp-thm}
that
$$
e^{2\nu} \Delta \tilde \vp + \tau _{\psi} \ge C \delta
> 0.
$$
Next, (\ref{dbar-eta}) and Lemma \ref{v-estimate} imply that
$\tilde \vp$ is comparable to $\vp$ on the support of $\dbar
\eta$, which lies in $V_{\sigma} (\gamma) := D_{\sigma} (\gamma) -
D_{\frac{\sigma}{ 2}}(\gamma)$.  We then have the estimate
\begin{eqnarray*}
\int _X |\dbar \eta |^2 e^{-2u_{\psi}}e^{-2\tilde \vp} & \le &
\frac{C}{\sigma ^2} \sum _{\gamma \in \Gamma} |s_{\gamma}|^2
e^{-2\vp (\gamma)} \int
_{V_{\sigma}(\gamma)} |\dbar \rho _{\gamma}|^2 e^{-2u_{\psi}}\\
& \le & \frac{C}{\sigma ^2} \sum _{\gamma \in \Gamma}
|s_{\gamma}|^2 e^{-2\vp (\gamma)} \int
_{D_{\sigma}(\gamma)} |\dbar \rho _{\gamma}|^2 e^{-2u_{\psi}}\\
& \le & C' \sum _{\gamma \in \Gamma} |s_{\gamma}|^2 e^{-2\vp
(\gamma)} \int
_{D_{\sigma}(\gamma)} e^{-2(\nu + u_{\psi})}\\
&<& +\infty,
\end{eqnarray*}
where the first inequality follows from Lemma \ref{local-fn} and the
last inequality follows from (\ref{ff-est}). Applying Theorem
\ref{o-dbar-thm}, we obtain a function $U \in L^2 (X,g,\tilde \vp)
\subset L^2(X,g,\vp)$ such that $\dbar U = \dbar \eta$. Moreover,
since $e^{-2\tilde \vp} \sim \frac{1}{|z-\gamma|^2}$ for $z$
sufficiently close to $\gamma$, we see that $U(\gamma)=0$ for all
$\gamma \in \Gamma$.  Thus the function
$$f := \eta - U \in \scb ^2 (X,g,\vp)$$
interpolates $(s_{\gamma})$, and the proof of Theorem
\ref{interp-thm} is complete.\qed

\subsection{Proof of Theorem
\ref{sampling-thm}}\label{sampling-section} Let $\hat \vp := \vp +
v_{r,\ve}$.  The main idea behind the proof of Theorem
\ref{sampling-thm} is the following sampling type lemma.

\begin{lem}\label{sampling-inequality}
Suppose the metric $e^{-2\psi}$ satisfies the differential
inequalities {\rm (\ref{diff-ineq})}.  For each $h \in \scb ^2
(X,g,\hat \vp )$,
\begin{eqnarray}\label{key-ineq-o}
\int _X |h|^2e^{-2\hat \vp} e^{2\nu} \Delta \hat \vp dA_g &\ge& 0.
\end{eqnarray}
\end{lem}

\pf Consider the function $S=|h|^2e^{-2\hat \vp }$. Then
\begin{eqnarray*}
\frac{\Delta S}{S}&=& \Delta \log  S + \frac{1}{S^2} |\di S|^2 =
\frac{1}{S^2} |\di S|^2 + \Delta \log |h|^2 - 2\Delta \hat \vp
\end{eqnarray*}
and thus
$$
e^{2\nu} \Delta S \ge - 2 S\left ( e^{2\nu}\Delta \hat \vp \right
) \ge  - 2 S\left ( e^{2\nu}\Delta \hat \vp \right ).
$$
We claim that
$$
\int _X e^{2\nu} \Delta S \ d\! A_g \le  0.
$$
To prove the claim, let $z_0 \in X$.  Take  $\lambda \in \cc
^{\infty} _0 ([0,1/2])$ such that $\lambda (t) \equiv 1$ for $0
\le t \le 1/4$, and put
$$\chi _a(r) := \lambda(r^2(1-a)).$$
Then
\begin{eqnarray*}
\int _X e^{-2(\psi - \nu)}\Delta S &=& \int _X e^{-2u_{\psi}}\Delta S\\
&=& \lim _{a \nearrow 1} \int _X e^{-2u_{\psi}} \chi _a \circ \rho
_{z_0}
\Delta S \\
&=& \lim _{a \nearrow 1}\int _X S \Delta \left (e^{-2u_{\psi}}
\cdot
(\chi _a \circ \rho _{z_0}) \right ) \\
&=& \lim _{a \nearrow 1} \int _X S \left ( \Delta (e^{-2u_{\psi}})
\chi _a \circ \rho _{z_0} + \di(e^{-2u_{\psi}}) \wedge \dbar (\chi
_a \circ
\rho _{z_0}) \right . \\
&& \qquad \qquad \left . + \dbar (e^{-2u_{\psi}}) \wedge \di (\chi
_a \circ
\rho _{z_0})+ e^{-2u_{\psi}} \Delta ( \chi _a \circ \rho _{z_0}) \right )\\
&=& \lim _{a \nearrow 1} \int _X S \left ( \Delta (e^{-2u_{\psi}})
\chi _a \circ \rho _{z_0} + e^{-2u_{\psi}} \Delta ( \chi _a \circ
\rho _{z_0}) \right ),
\end{eqnarray*}
where the third equality follows from Stokes' Theorem. Now,
\begin{eqnarray*}
&& \lim _{a \nearrow 1} \int _X Se^{-2u_{\psi}} \Delta ( \chi _a
\circ \rho _{z_0}) \\
&=& \lim _{a \nearrow 1} \int _X S e^{-2u_{\psi}} \left (\chi _a
'' (\rho _{z_0})|\di \rho _{z_0}|^2 +
\chi _a ' (\rho _{z_0}) \Delta \rho _{z_0} \right )\\
&=& \lim _{a \nearrow 1} \int _X S e^{-2u_{\psi}} \left (\chi _a
'' (\rho _{z_0})
+ \frac{\chi _a ' (\rho _{z_0})}{\rho _{z_0}} \right )|\di \rho _{z_0}|^2\\
&=& \lim _{a \nearrow 1} \int _X |h|^2 e^{-2\hat \vp} \left (\chi
_a '' (\rho _{z_0}) + \frac{\chi _a ' (\rho _{z_0})}{\rho _{z_0}}
\right ) e^{2\nu}|\di \rho _{z_0}|^2 dA_g =0,
\end{eqnarray*}
where the last equality follows from (\ref{ff-est}) and the
definition of $\chi _a$.  It follows that
 \begin{eqnarray*}
\int _X e^{-2(\psi - \nu)}\Delta S&=& \int _X S\Delta
e^{-2u_{\psi}}.
\end{eqnarray*}
Since
$$
\Delta e^{-2u_{\psi}} = 2 e^{-2u_{\psi}}(2|\di u_{\psi}|^2 -
\Delta u_{\psi}) = 2 e^{-2\psi}e^{2\nu} \left (2|\di u_{\psi}|^2 -
\Delta u_{\psi}\right ) ,
$$
the lemma now follows from (\ref{diff-ineq}). \qed

\medskip

\noi {\it Conclusion of the proof of Theorem \ref{sampling-thm}.}
Let $h \in \scb ^2 (X,g,\vp)$.  By Lemma \ref{v-estimates} we
calculate that
\begin{eqnarray*}
e^{2\nu(z)} \Delta \hat \vp(z) &=& e^{2\nu(z)} \Delta \vp(z) +
e^{2\nu (z)} \Delta v_{r,\ve}(z)\\
&=&  e^{2\nu(z)} \Delta \vp(z) \left ( 1 - t \sum _{\gamma \in
\Gamma} \frac{1}{2\ve ^2 } \int _{D_{\ve}(\gamma)}\frac{\xi _r
(\zeta,z)}{e^{2\nu(z)} \Delta \vp(z)}dA_{E,\gamma}(\zeta) \right . \\
&& \qquad \qquad  \qquad  \qquad  \left . + t \sum _{\gamma \in
\Gamma} e^{2\psi} \frac{1}{\ve ^2} \frac{e^{2\psi(z)}|\di \rho
_{\gamma}(z) |^2}{e^{2\nu(z)} \Delta \vp(z)} {\bf 1}_{D_
{\ve}(\gamma)}(z) \right ).
\end{eqnarray*}
Applying the hypotheses $D^-_f(\Gamma) >1$ and (\ref{ff-est}), we
see therefore that, for $t$ sufficiently close to $1$, there exist
$r,\delta,C
>0$ such that
\begin{eqnarray}\label{r-dens}
e^{2\nu} \Delta \hat \vp \le - t e^{2\nu} \Delta \vp\left ( \delta -
C \sum _{\gamma \in \Gamma}  e^{2\psi} \frac{2}{\ve ^2} {\bf 1}_{D_
{\ve}(\gamma)} \right ).
\end{eqnarray}
We then apply Lemma \ref{sampling-inequality} to get
\begin{eqnarray*}
\int _X |h|^2 e^{-2\vp} dA_g &\le &\int _X |h|^2 e^{-2\hat \vp}
dA_g \\
&\le& C \int _X e^{2\nu}\Delta (\vp )|h|^2 dA_g \\
&\le& C' \sum _{\gamma \in \Gamma} \frac{2}{\ve ^2} \int _{D_{\ve}
(\gamma )} e^{2\nu}\Delta (\vp)|h|^2e^{-2
\hat \vp} dA_g \\
&\le & C'' \sum _{\gamma \in \Gamma} \frac{2}{\ve ^2} \int
_{D_{\ve} (\gamma )}|h|^2e^{-2\hat \vp} dA_g \\
&\le& C''' \sum _{\gamma \in \Gamma} \frac{2}{\ve ^{2+2t}} \int
_{D_{\ve} (\gamma )}|h|^2e^{-2\vp} dA_g,
\end{eqnarray*}
where the first  inequality follows from Lemma \ref{v-estimates},
the third inequality follows from integration of (\ref{r-dens})
together with Lemma \ref{sampling-inequality} and the last
inequality follows from Lemmas \ref{v-estimates} and \ref{min}.
Now,
\begin{eqnarray*}
\int _{D_{\ve} (\gamma )}\!\!\!\!\!\!\!\! |h|^2e^{-2\vp} dA_g &=&
\int _{D_{\ve} (\gamma )}\!\!\!\!\!\!\!\!
|he^{-F_{\gamma}}|^2e^{-2\vp+ 2{\rm Re}F_{\gamma}} dA_g \\
&\le & C e^{-2\vp(\gamma)} \int
_{D_{\ve} (\gamma )}\!\!\!\!\!\!\!\! |he^{-F_{\gamma}}|^2 dA_g \\
&\le& C' A_g (D_{\ve} (\gamma) ) e^{-2\vp (\gamma)}\!\! \left (
|h(\gamma)|^2 \! + \! \ve ^2 \sup _{D_{\ve} (\gamma)} \frac{\left |
(he^{-F_{\gamma}})' \right
|^2}{|\di \rho _{\gamma}|^2}\right )\\
&\le & C' A_g (D_{\ve} (\gamma) )
e^{-2\vp (\gamma)} \\
&&\qquad \qquad \times \left ( |h(\gamma)|^2 + \ve ^2C_{\ve,\sigma}
\int _{D_{\sigma}
(\gamma)}\!\!\!\!\!\!\!\!\!\! |he^{-F_{\gamma}}|^2 dA_g \right )\\
&\le & C' A_g (D_{\sigma} (\gamma) )
e^{-2\vp (\gamma)} |h(\gamma)|^2 \\
&&\qquad \qquad + \ve ^2 C'' A_g(D_{\ve}(\gamma)) \int _{D_{\sigma}
(\gamma)}\!\!\!\!\!\!\!\!\!\! |h|^2 e^{-2\vp} dA_g, \\
\end{eqnarray*}
where the first and last inequalities follow from Lemma
\ref{local-fn}, the second inequality follows from Taylor's
theorem, and the third inequality from the Cauchy estimate
(\ref{der-est}).

Next, since $X$ is fundamentally finite and $e^{-2\psi} \le
e^{-2\nu}$, we see that
$$
A_g(D_{\ve}(\gamma)) \le \int _{D_{\ve}(\gamma)} \!\!\!\!\!\!\!\!
e^{-2\nu}\le C \int _{D_{\ve}(\gamma)} \!\!\!\!\!\!\!\! |\di \rho
_{\gamma}|^2 = \pi C \ve ^2
$$
for all sufficiently small $\ve$ and some $C$ independent of
$\gamma$, where the last equality follows from (\ref{exact-area}).
We thus obtain
\begin{eqnarray*}
&& \int _X \!\! |h|^2 e^{-2\vp}dA_g \\
&& \le \sum _{\gamma \in \Gamma} \left ( \frac{C_1}{\ve ^{2+2t}}
|h(\gamma)|^2e^{-2\vp(\gamma)} A_g(D_{\sigma}(\gamma)) + C_2 \ve
^{2-2t} \!\!\int _{D_{\sigma} (\gamma)}\!\!\!\!\!\!\!\!\!\!
|h|^2 e^{-2\vp}dA_g \right ) \\
&& \le \sum _{\gamma \in \Gamma}  \left ( \frac{C_1}{\ve ^{2+2t}}
|h(\gamma)|^2e^{-2\vp(\gamma)} A_g(D_{\sigma}(\gamma))\right ) + C_2
\ve ^{2-2t} \int _X |h|^2 e^{-2\vp}dA_g.
\end{eqnarray*}
By taking $\ve$ sufficiently small, we obtain the left hand side of
(\ref{sampling-ineq}).  For the right hand side of
(\ref{sampling-ineq}), we argue as follows.
\begin{eqnarray*}
\sum _{\gamma \in \Gamma} |h(\gamma)|^2e^{-2\vp(\gamma)}
A_g(D_{\sigma}(\gamma)) &=& \sum _{\gamma \in
\Gamma}|h(\gamma)e^{-F_{\gamma}(\gamma)}|^2
e^{-2\vp(\gamma)} A_g(D_{\sigma}(\gamma)) \\
&\le& C \sigma ^2 \sum _{\gamma \in \Gamma} e^{-2\vp(\gamma)}
\int _{D_{\sigma}(\gamma)}\!\!\!\!\!\!\!\!\!\! |he^{-F_{\gamma}}|^2 dA_g \\
&\le & C' \sum _{\gamma \in \Gamma} \int _{D_{\sigma}
(\gamma)}\!\!\!\!\!\!\!\!\!\! |h|^2 e^{-2\vp}dA_g \\
&\le& C'' \int _X |h|^2 e^{-2\vp}dA_g,
\end{eqnarray*}
where the first inequality follows from (\ref{est}), the second
from Lemma \ref{local-fn} and the third from the definition of the
separation constant.  The proof of Theorem \ref{sampling-thm} is
thus complete.\qed

\section{Examples}\label{examples-section}

\subsection{The Euclidean plane}

In this paragraph, we consider the case of the Euclidean complex
plane $(X,g)=(\C, |dz|^2)$.  The generalized Bergman space in this
situation is
$$
\scb \scf ^2 = \left \{ h \in \sco (\C)\ ;\ ||h||^2_{\vp} := \int
_\C |h|^2 e^{-2\vp}dm < +\infty \right \},
$$
where $dm$ is Lebesgue measure in the plane, and
$$
\fb \ff ^2 = \left \{ (s_{\gamma}) \subset \C \ ;\
||(s_{\gamma})||^2_{\vp} := \sum _{\gamma \in \Gamma} |s_{\gamma}|^2
e^{-2\vp(\gamma)} < + \infty \right \}.
$$
The space $\scb \scf ^2$ is  sometimes called {\it generalized
Bargmann-Fock space}. When $\vp (z) = |z|^2/2$ we obtain the
classical Bargmann-Fock space.

The plane is the main example of a parabolic Riemann surface.  The
Evans kernel in $\C$ is unique and is given by $E(z,\zeta) = \log
|z-\zeta|.$ Thus $\rho _z(\zeta) = |z-\zeta|$ and the disks
$D_{\sigma}(z)$ are simply the Euclidean disks $|z-\zeta| < \sigma$.
A simple calculation shows that
$$
|d\rho _z(\zeta)|^2 = 4 |\di \rho _z(\zeta)|^2 = 1,
$$
and thus the fundamental metric is just a multiple of the Euclidean
metric.

The upper and lower densities are given by
$$
D^+_f (\Gamma ) = \limsup _{r \to \infty} \sup _{z \in \C} \sum
_{\Gamma \cap D_r(z) } \frac{f(|z-\gamma|)}{4 \Delta \vp \int _0
^r tf(t)dt}
$$
and
$$
D^-_f (\Gamma)  = \liminf _{r \to \infty} \inf _{z \in \C} \sum
_{\Gamma \cap D_r(z) } \frac{f(|z-\gamma|)} {4 \Delta \vp \int _0
^r tf(t)dt}.
$$

If we choose as our locally integrable function $f$ the
constant function, we recover the results of \cite{quimbo}.
However, by making other choices, we can get other sufficient
conditions that, although not necessary, might be of use in some
applications.

For the sake of simplicity, we will consider in the following examples only
the classical Bargmann-Fock space.

\begin{exa}
\begin{enumerate}
\item[(i)]
Let $f(t)=e^{-t}$. Then $\Gamma$ is interpolating if
$$
 \sup _{z \in \C} \sum
_{\Gamma \cap D_r(z) } e^{-|z-\gamma|}<2
$$
and sampling if
$$
\inf _{z \in \C} \sum
_{\Gamma \cap D_r(z)} e^{-|z-\gamma|}>2.
$$
Integration by parts, together with a standard argument shows that
$\Gamma$ is interpolating if
$$
\sup_{z\in\C} \int_0^{\infty} \# (\Gamma
\cap D_s(z))\frac{ds}{e^s}<2
$$
and sampling if
$$
\inf_{z\in\C} \int_0^{\infty} \# (\Gamma
\cap D_s(z))\frac{ds}{e^s}>2.
$$

\item[(ii)] Let $f_{a} :=
{\bf 1}_{[0,a]}$. We then obtain:

 If  $a > 1/\sqrt{2}$ and every  disk of radius $a$ contains at most
one member of $\Gamma$, then $\Gamma$ is interpolating.

If $a < 1/\sqrt{2}$ and every disk of radius
$a$ contains at least one member of $\Gamma$, then $\Gamma$ is
sampling.
\end{enumerate}
\end{exa}

\subsection{The disk}

In this paragraph we consider the case of the Poincar\' e unit disk
$(X, g)=(\D, \frac{|dz|^2}{(1-|z|^2)^2})$.  The disk is the main
example of a regular hyperbolic Riemann surface.  Its Green's
function is
$$E(z,\zeta) = \log |\phi _z (\zeta)|, \qquad {\rm where} \quad
\phi _z(\zeta) = \frac{z-\zeta}{1-\bar z \zeta}$$
is the standard involution.  Thus $\rho _z(\zeta) = |\phi _z (\zeta)|$ and
the disks $D_{\sigma} (z)$ are the well-known pseudo-hyperbolic disks; they
are geometrically Euclidean disks, but their Euclidean centers and radii are
different.

Standard calculations show that
$$
|d \rho _z(\zeta)|^2 = \left | \frac{1-|z|^2}{(1-\bar z \zeta)^2} \right |^2
= \frac{(1-\rho _z(\zeta)^2)^2}{(1-|\zeta|^2)^2},
$$
so we have $\nu(\zeta)=\log2 +\log(1-|\zeta|^2)$, and it is clear
that (\ref{ff-est}) holds.

As suggested by the proof of Proposition \ref{some-u}, we take
$$
u_{\psi} (z) = - \frac{1}{2} \log (1-|z|^2),
$$
and thus we have
$$
e^{2\nu(z)}(\Delta u_{\psi}(z) - 2|\di u _{\psi}(z)|^2) =
\frac{1}{2}(1-|z|^2)\ge 0 \qquad {\rm and} \qquad \tau _{\psi}(z)
=\frac{1}{2(1-|z|^2)}.
$$
We also have
$$
A_g(D_{\sigma}(\gamma)) = C_{\sigma} (1-|\gamma|^2).
$$
Thus our Hilbert spaces are
$$
\cb _{\vp} ^2 := \left \{ h \in \sco (\D)\ ;\ \int _{\D} |h|^2
e^{-2\vp} \frac{d m} {(1-|z|^2)} < +\infty \right \}
$$
and
$$
b^2_{\vp} := \left \{ (s_{\gamma})\ ;\ \sum _{\gamma \in \Gamma}
|s_{\gamma}|^2 e^{-2\vp(\gamma)}(1-|\gamma|^2) < +\infty \right
\}.
$$
The densities are given by
\begin{eqnarray*}
D_f^ + (\Gamma) &=& \limsup _{r\to1} \sup_{z\in\D}\sum_ {\rho_z
(\gamma)<r}\frac{f(\rho_z(\gamma)) (1-\rho_z(\gamma)^2)^2} {4\left
( (1-|z|^2)^2\Delta \vp (z) + \frac{1}{2}(1-|z|^2) \right
)\int_0^r
tf(t)dt},\\
&& {\rm and}  \\
D_f^- (\Gamma) &=& \liminf _{r\to1} \inf _{z\in\D}
\sum_{\rho_z(\gamma)<r} \frac{f(\rho_z(\gamma)) (1-\rho_z(\gamma)
^2)^2}{4(1-|z|^2)^2 \Delta \vp(z) \int_0^r tf(t)dt}.
\end{eqnarray*}

If we take
$$
f(t)=\frac{-\log t}{(1-t^2)^2}{\bf 1}_{\left [ \tfrac{1}{2}\ ,\ 1
\right )},
$$
Theorems \ref{interp-thm} and \ref{sampling-thm} recover the
results from \cite{quimbo}.

Again for the sake of illustration we will consider below only
the classical unweighted Bergman space, which is obtained by setting
 $\vp = -\frac{1}{2}\log (1-|z|^2)$.

\begin{exa}
\begin{enumerate}
\item [(i)]  Letting $f=1$, we see that $\Gamma$ is
interpolating if
$$\sup_{z\in\D} \sum(1- \rho_z (\gamma)^2)^2<1$$
and sampling if
$$\inf _{z\in\D} \sum (1-\rho_z(\gamma)^2 )^2>1.$$

\item[(ii)] Letting $f(t)=(1-t^2)^{-2}$, we see that
$\Gamma$ is interpolating if
$$
\limsup _{r\to1} \sup_{z\in\D} \frac{\# ( \Gamma \cap
D_r(z))}{A_{\rm hyp}(D_r(z))}<1
$$
and sampling if
$$\liminf _{r\to1} \inf_{z\in\D}
\frac{\#(\Gamma \cap D_r(z))}{A_{\rm hyp}(D_r(z))}>1,
$$
where
$$
A_{\rm hyp}(D_r(z)) = \frac{1}{2\pi}
\int _{D_r(z)} \frac{dm(z)}{(1-|z|^2)^2}
$$
denotes hyperbolic area of $D_r(z)$.

\item[(iii)] Let $f_{a} :=
{\bf 1}_{[0,a]}$. We then obtain:

 If $\delta > \frac{1}{\sqrt{2}}$ and $\Gamma$ has
at most one point in every disk of radius $\delta$, then $\Gamma$
is interpolating.

 If $\delta < \frac{1}
{\sqrt{2}}$ and every disk of radius $\delta$ contains at
least one member of $\Gamma$, then $\Gamma$ is sampling.
\end{enumerate}
\end{exa}

\end{document}